\documentclass[twoside,11pt]{article}
\usepackage{amstex,amssymb}
\usepackage{theorem}
\usepackage{times}
\usepackage{a4}
{\theorembodyfont{\slshape}
\newtheorem{theorem}{Theorem}[section]
\newtheorem{lemma}[theorem]{Lemma}
\newtheorem{proposition}[theorem]{Proposition}
\newtheorem{corollary}[theorem]{Corollary}
\newtheorem{conjecture}[theorem]{Conjecture}
}
{\theorembodyfont{\upshape}
\newtheorem{definition}[theorem]{Definition}
\newtheorem{example}[theorem]{Example}

\newtheorem{remark}[theorem]{Remark}

\newtheorem{notit}[theorem]{{}}

}

\newcommand{\prf}{{\kursiv Proof}. }

\newcommand{\prfofthetheorem}{{\kursiv Proof of the theorem}. }

\newcommand{\qed}{\hspace*{\fill}$\Box$\vspace{2ex}}

\newcommand{\beeq}[1]{\begin{eqnarray}\label{#1}}
\newcommand{\eneq}{\end{eqnarray}}

\newcommand{\kursiv}{\sl}

\newcommand{\ka}{{\cal A}}

\newcommand{\ke}{{\cal E}}

\newcommand{\kh}{{\cal H}}
\newcommand{\ki}{{\cal I}}

\newcommand{\kn}{{\cal N}}
\newcommand{\ko}{{\cal O}}

\newcommand{\kt}{{\cal T}}

\newcommand{\IA}{{\mathbb A}}

\newcommand{\IC}{{\mathbb C}}

\newcommand{\IH}{{\mathbb H}}

\newcommand{\IP}{{\mathbb P}} 
\newcommand{\IQ}{{\mathbb Q}} 
\newcommand{\IR}{{\mathbb R}}
 
\newcommand{\IZ}{{\mathbb Z}}

\newcommand{\gotha}{{\mathfrak a}}
\newcommand{\gothb}{{\mathfrak b}}
\newcommand{\gothc}{{\mathfrak c}}
\newcommand{\gothd}{{\mathfrak d}}
\newcommand{\gothe}{{\mathfrak e}}
\newcommand{\gothf}{{\mathfrak f}}
\newcommand{\gothg}{{\mathfrak g}}
\newcommand{\gothh}{{\mathfrak h}}

\newcommand{\gothm}{{\mathfrak m}}

\newcommand{\gothq}{{\mathfrak q}}

\newcommand{\gothy}{{\mathfrak y}}

\newcommand{\gothA}{{\mathfrak A}}
\newcommand{\gothB}{{\mathfrak B}}
\newcommand{\gothC}{{\mathfrak C}}

\newcommand{\gothL}{{\mathfrak L}}

\newcommand{\gothS}{{\mathfrak S}}
\newcommand{\gothX}{{\mathfrak X}}

\newcommand{\dual}{\makebox[0mm]{}^{{\scriptstyle\vee}}}

\newcommand{\isom}{\cong}

\newcommand{\tensor}{\otimes}

\newcommand{\Hilb}{{\rm Hilb}}

\newcommand{\Supp}{{\rm Supp}}

\newcommand{\Spec}{{\rm Spec}}

\newcommand{\End}{{\rm End}}

\newcommand{\cc}{\gothc}
\newcommand{\dd}{\gothd}
\newcommand{\ee}{\gothe}

\newcommand{\qq}{\gothq}
\newcommand{\LL}{\gothL}
\newcommand{\ch}{{\gothc\gothh}}
\newcommand{\ad}{{\rm ad\,}}
\newcommand{\vacuum}{{\bf 1}}

\renewcommand{\phi}{\varphi}
\newcommand{\blowup}[2]{{\rm Bl}_{#1}{(#2})}

\newcommand{\id}{{\rm id}}

\newcommand{\rk}{{\rm rk}}

\newcommand{\Coh}{{\rm Coh}}

\newcommand{\length}{{\rm length}}

\newcommand{\pr}{{\rm pr}}

\newcommand{\lra}{\longrightarrow}
\newcommand{\lla}{\longleftarrow}
\newcommand{\ses}[3]{0\rightarrow#1\rightarrow#2
   \rightarrow#3\rightarrow0}

\newcommand{\verylongrightarrow}[1]{\raisebox{0mm}[0.9ex][0ex]{\hbox to #1{\rightarrowfill}}}
\newcommand{\verylongleftarrow}[1]{\raisebox{0mm}[0.9ex][0ex]{\hbox to #1{\leftarrowfill}}}
\newcommand{\rpfeil}[2]{\stackrel{#2}{\verylongrightarrow{#1}}}

\newcommand{\rationalmap}{{-{}-{}\!\to}}

\begin{document}
\title{Chern Classes of Tautological Sheaves on \linebreak Hilbert Schemes of Points on Surfaces}
\author{Manfred Lehn}
\date{}
\maketitle

\begin{abstract}
\noindent
We give an algorithmic description of the action of the Chern classes of tautological bundles on the cohomology of Hilbert schemes of points on a smooth surface within the framework of Nakajima's oscillator algebra.
This leads to an identification of the cohomology ring of $\Hilb^n(\IA^2)$
with a ring of explicitly given differential operators on a Fock space. We
end with the computation of the top Segre classes of tautological bundles associated to line bundles on $\Hilb^n$ up to $n=7$, extending computations of Severi, LeBarz, Tikhomirov and Troshina and give a conjecture for the generating series.
\end{abstract}

\section*{Introduction}

Hilbert schemes $X^{[n]}$ of $n$-tuples of points on a complex 
projective manifold $X$ 
are natural compactifications of the configuration spaces of unordered
distinct $n$-tuples of points on $X$. Their geometry is determined by the
geometry of $X$ itself and the geometry of the `punctual' Hilbert schemes
of all zero-dimensional subschemes in affine space that are supported at the
origin. Thus one is naturally led to the following problem:

{\em Determine
explicitly the geometric or topological invariants of the Hilbert schemes 
$X^{[n]}$ such as the Betti numbers, the Hodge numbers, the Chern numbers, the 
cohomology ring, from the corresponding data of the manifold $X$ itself.}

This problem is most attractive when $X$ is a surface, since then the Hilbert schemes are themselves irreducible projective manifolds,
by a result of  Fogarty  \cite{Fogarty},
whereas for higher dimensional varieties the Hilbert schemes are in general
neither irreducible nor smooth nor pure of expected dimension. 

The answer to the problem above for the Betti numbers was given for $\IP^2$ and rational ruled surfaces by Ellingsrud and Str\o{}mme \cite{EllStr1} and for general surfaces by G\"ottsche in \cite{Goettsche}. The answer turns out to be 
particularly beautiful (cf.\ Theorem \ref{GoettschesFormel} below).
The problem for the Hodge numbers was solved by  S\"orgel  and  G\"ottsche
\cite{GoettscheSoergel}. For a different approach to both results see 
\cite{Cheah}. A partial answer for the Chern classes will be given in a forthcoming paper by Ellingsrud, G\"ottsche and the author \cite{EllingsrudGoettscheLehn}.

The question for the ring structure of the cohomology is more difficult. In general, $X^{[2]}$ is the quotient of the blow-up of $X\times X$ along the diagonal by the canonical involution that exchanges the factors. Thus the case
of interest is $H^*(X^{[n]})$, $n\geq3$. The ring structure for
$H^*(X^{[3]})$, $X$ smooth projective of arbitrary dimension, was found
by Fantechi and  G\"ottsche \cite{FantechiGoettsche}.
In another direction, Ellingsrud and Str\o{}mme \cite{EllStr2} gave generators 
for $H^*((\IP^2)^{[n]},\IZ)$, $n$ arbitrary, and an implicit description of the relations.

Vafa and  Witten \cite{VafaWitten} remarked that G\"ottsche's Formula for the 
Betti numbers is identical with the Poincar\'e series of a Fock space modelled 
on the cohomology of $X$. Nakajima \cite{Nakajima} succeeded in giving a 
geometric construction of such a Fock space structure on the cohomology of the 
Hilbert schemes, leading to a natural `explanation' of G\"ottsche's result. 
Similar results have been announced by Grojnowski \cite{Grojnowski}.

Following the presentation of  Grojnowski, this can be made more precise
as follows: sending a pair $(\xi',\xi'')$ of subschemes
of length $n'$ and $n''$, respectively, and of disjoint support to their
union $\xi'\cup\xi''$ defines a rational map 
$$m:X^{[n']}\times X^{[n'']}\rationalmap X^{[n'+n'']}.$$
This map induces linear maps on the rational cohomology
$$m_{n',n''}:H^*(X^{[n']};\IQ)\tensor H^*(X^{[n'']};\IQ)\lra H^*(X^{[n'+n'']};\IQ)$$
and
$$m^{n',n''}:H^*(X^{[n'+n'']};\IQ)\lra H^*(X^{[n']};\IQ)\tensor H^*(X^{[n'']};\IQ).$$
If we let $\IH:=\oplus_nH^*(X^{[n]};\IQ)$, then these maps define a multiplication
and a comultiplication
$$m_*:\IH\tensor\IH\lra \IH,\qquad m^*:\IH\lra \IH\tensor\IH,$$
which make $\IH$ a commutative and cocommutative bigraded Hopf algebra. The 
result of Nakajima and Grojnowski says that this Hopf algebra is isomorphic 
to the graded symmetric algebra of the vector space $H^*(X;\IQ)\tensor t\IQ[t]$. 

More explicitly, Nakajima constructed linear maps
\footnote{Our presentation differs in notations and conventions slightly from
Nakajima's.}
$$\qq_n:H^*(X;\IQ)\lra\End_\IQ(\IH),\quad n\in\IZ,$$
and proved that they satisfy the `oscillator' or `Heisenberg' relations
$$[\qq_n(\alpha),\qq_m(\beta)]= n\cdot 
\delta_{n+m}\cdot\int_X\alpha\beta\cdot\id_{\IH}.$$
Here the commutator is to be taken in a graded sense.

The multiplication and the comultiplication of $\IH$ are not obviously related
to the quite different ring structure of $\IH$, which is given by the usual
cup product on each direct summand $H^*(X^{[n]};\IQ)$. (Strictly speaking, $\IH$
contains a countable number of idempotents $1_{X^{[n]}}\in H^0(X^{[n]};\IQ)$ but
not a unit unless we pass to some completion).

This paper attempts to relate
the Hopf algebra structure and the cup product structure. More precisely:

Let $F$ be locally free sheaf of rank $r$ on $X$. Attaching to a point
$\xi\in X^{[n]}$, i.e.\ a zero-dimensional subscheme $\xi\subset X$, the
$\IC$-vector space $F\tensor \ko_\xi$ defines a locally free sheaf $F^{[n]}$
of rank $rn$ on $X^{[n]}$. The Chern classes of all sheaves on $X^{[n]}$
of this type generate a subalgebra $\ka\subset\IH$. We will describe a
purely algebraic
algorithm to determine the action of $\ka$ on $\IH$ in terms of the $\IQ$-basis
of $\IH$ provided by Nakajima's results.
We collect the Chern classes of all sheaves $F^{[n]}$ for
a given sheaf $F$ into operators
$$\ch(F):\IH\to \IH,\qquad\cc(F):\IH\to\IH$$
and geometrically compute the commutators of these operators with the
oscillator operators defined by Nakajima. 

A central r\^ole is played by the operator $\dd:=\cc_1(\ko_X)$, which ---
up to a factor $(-1/2)$ --- can also be interpreted as the intersection with
the `boundaries' of the Hilbert schemes, i.e.\ the divisors $\partial X^{[n]}\subset X^{[n]}$ of all tuples $\xi$ which have a multiple point somewhere. The 
derivative of any operator $\gothf\in\End(\IH)$ is defined by $\gothf':=
[\dd,\gothf]$. Our main technical result then says that for $n>0$
\beeq{erstesErgebnis}
\qq_n'(\alpha)=\frac{n}{2}\sum_{\nu}\qq_\nu\qq_{n-\nu}\delta(\alpha)+\binom{n}{2}\qq_n(K\alpha),
\eneq
where $\delta:H^*(X;\IQ)\to H^*(X;\IQ)\tensor H^*(X;\IQ)$ is the map induced 
by the diagonal embedding and $K$ is the canonical class of $X$. An immediate 
algebraic consequence of this relation is
\beeq{zweitesErgebnis}
[\qq_n'(\alpha),\qq_m(\beta)]=-nm\cdot\qq_{n+m}(\alpha\beta)
\eneq
for $n,m>0$. By induction one concludes that the operators $\qq_1$ and $\dd$
suffice to generate all $\qq_n$, $n\geq 1$. 

The commutator of the Chern character operator $\ch(F)$ with the 
standard operator $\qq_1$ can be expressed in terms of higher derivatives of 
$\qq_1$:
\beeq{drittesErgebnis}
[\ch(F),\qq_1(\alpha)]=\sum_{n\geq0}\frac{1}{n!}\qq_1^{(n)}(ch(F)\alpha).
\eneq
Equations (\ref{erstesErgebnis}), (\ref{zweitesErgebnis}) and (\ref{drittesErgebnis}) together give a complete description of the action of 
$\ka$ on $\IH$. Here are some applications:

1. We prove the following formula conjectured by G\"ottsche:
If $L$ is a line bundle on $X$ then
$$\sum_{n\geq0}c(L^{[n]})z^n=\exp\left(\sum_{m\geq1}\frac{(-1)^{m-1}}{m}\qq_m(c(L))z^m\right).$$

2. We give a general algebraic solution to Donaldson's question for the
integral $N_n$ of the top Segre class of the bundles $L^{[n]}$ associated to a line bundle $L$ for any $n$ and explicitly compute $N_n$ for $n\leq7$. From
an analysis of this computational material we derive a conjecture for the
generating function for all $N_n$.

3. We identify the Chow ring of the Hilbert scheme of the affine plane
with an algebra of explicitly given differential operators on the polynomial
ring $\IQ[q_1,q_2,\ldots]$ of countably many variables.

This paper is organised as follows: In Section \ref{Preliminaries} we recall
the basic geometric notions used in the later parts. Section \ref{TheStructure} 
provides an introduction to Nakajima's results. Section \ref{TheDerivative} 
contains the core of this paper: we first define Virasoro operators $\LL_n$ in
analogy to the standard construction and show how these arise geometrically.
We then introduce the operator $\dd$ and compute the derivative of $\qq_n$.
Finally, in Section \ref{RingStructure} we apply these results to compute
the action of the Chern classes of tautological bundles.
 
Discussions with A.\ King were
important to me in clarifying and understanding the picture that Nakajima
draws in his very inspiring article. I am very grateful to G.\ Ellingsrud
for all the things I learned from his talks and conversations with him about
Hilbert schemes. To some extend the results in this article are a reflection on
an induction method entirely due to him. I thank W.\ Nahm for pointing out
a missing factor in Theorem \ref{Lqcommutator} and D.\ Zagier for a very instructive correspondence on power series.

Most of the research for this paper
was carried out during my stay at the SFB 343 of the University of Bielefeld.
I owe special thanks to S.\ Bauer for his continuous encouragement, interest and support.


\section{Preliminaries}\label{Preliminaries}

In this section we introduce the basic notations that will be used throughout
the paper and collect some results from the literature without proof.
All varieties and schemes are of finite type over the complex numbers. $X$
will always denote a smooth irreducible projective surface.
If $f:S\to S'$ is a morphism of schemes, I will write $f_X:=(f\times\id_X):S\times X\to S'\times X$.

\subsection{Hilbert schemes of points}

For any smooth projective surface $X$ let $S^nX$ denote the symmetric product,
i.e.\ the quotient of $X^n$ by the action of the symmetric group $\gothS_n$,
and let $X^{[n]}$ be the Hilbert scheme of zero-dimensional closed subschemes of
length $n$. By a result of Grothendieck \cite{Bourbaki221} $X^{[n]}$ is again
a projective scheme. There is a natural morphism $\rho:X^{[n]}\to S^nX$, the
{\sl Hilbert-Chow} morphism, which maps a point $[\xi]\in X^{[n]}$ to the 
cycle $\sum_x\ell(\ko_{\xi,x})\cdot x$ (cf.\ Iversen \cite{Iversen}). 

The basic geo\-metry of the Hilbert schemes of points on surfaces is 
governed by two theorems due to Fogarty \cite{Fogarty} and Brian\c{c}on 
\cite{Briancon}.  

\begin{theorem}[Fogarty]--- $X^{[n]}$ is a $2n$-dimensional irreducible smooth variety.
\end{theorem}

\begin{remark}--- If $C$ is a curve, its Hilbert scheme $C^{[n]}$ is smooth
and the map $\rho:C^{[n]}\to S^nC$ is an isomorphism. Computing the dimension of the
tangent spaces one can show that $Y^{[3]}$ is smooth for a smooth variety $Y$ of any dimension. On the other hand, $Y^{[n]}$ is singular if $\dim(Y)>2$ and 
$n>3$.\qed
\end{remark}

Fix a point $p\in X$ and let $X^{[n]}_p\subset X^{[n]}$ 
denote the closed subset of all subschemes $\xi\subset X$ with $\Supp(\xi)=\{p\}$ (with the reduced induced subscheme structure). This is
indeed a closed subset, as it is the fibre $\rho^{-1}(np)$ of the
Hilbert-Chow morphism over the point $np\in S^nX$.

Let $(\ko,\gothm)$ denote
the local ring of $X$ at $p$. Since any point $\xi\in X^{[n]}_p$ may be 
considered as a subscheme of $\Spec(\ko/\gothm^n)$, and since
$\ko/\gothm^n\isom\IC[x,y]/(x,y)^n$, all schemes $X^{[n]}_p$ --- for varying $X$
and $p$ --- are (non-canonically) isomorphic. Clearly, $X^{[1]}_p=\{p\}$ and 
$X^{[2]}_p=\IP(T_{p}X\dual)$, moreover it is not too difficult to see that
$X^{[3]}_p$ is isomorphic to the projective cone over the twisted cubic
$C_3\subset \IP^3$, the vertex of the cone corresponding to the subscheme
$\Spec(\ko/\gothm^2)$. It is not accidental that in these examples the
dimension of $X^{[n]}_p$ increases by one in each step:

\begin{theorem}[Brian\c{c}on]--- For all $n\geq1$, $X^{[n]}_p$ is an irreducible variety of dimension $n-1$.\qed
\end{theorem}

For a proof see \cite{Briancon}. A new proof with a more geometric and
conceptual argument was recently given by Ellingsrud and Str\o{}mme
\cite{EllStrCoefficient}.

Brian\c{c}on's Theorem emphasises the importance of curvilinear schemes: recall
that a zero-dimensional subscheme $\xi\subset X$ is called {\em curvilinear} 
at $x\in X$, if $\xi_x$ is contained in some smooth curve $C\subset X$. Equivalently, one might say that $\ko_{\xi,x}$ is isomorphic to the 
$\IC$-algebra $\IC[z]/(z^\ell)$, where $\ell=\ell(\xi_x)$.
Hence $\xi$ is curvilinear at $x$ if $\xi_x$ is either empty, a reduced point, or if 
$\dim T_x\xi=1$. From this criterion it is clear, that in any flat family of zero-dimensional subschemes the points in the base space which correspond to curvilinear subschemes form an open subset.

In particular, we may consider the open subset $X^{[n]}_{p,curv}\subset 
X^{[n]}_p$. This set has a very nice structure:

\begin{lemma}\label{curvilinearones}--- If $n\geq 2$, then the morphism
$$t:X^{[n]}_{p,curv}\lra \IP(T_pX\dual), [\xi]\mapsto [T_p\xi]$$
is a bundle morphism with affine fibres $\IA^{n-2}$. In particular, $X^{[n]}_{p,curv}$ is an irreducible smooth variety of dimension $n-1$.
\end{lemma}

\prf Let $x,y\in\ko_{X,p}$ be local coordinates and consider the open subset
$U=\{\langle y+\alpha_1x\rangle|\alpha_1\in \IC\}\subset \IP(T_pX\dual)$. Then there is an isomorphism $\IA^{n-1}\to t^{-1}(U)$ sending the $(n-1)$-tuple 
$(\alpha_1,\ldots,\alpha_{n-1})$ to the subsheaf corresponding
to the ideal $(y+\alpha_1x+\ldots+\alpha_{n-1}x^{n-1})+\ki_{p}^{n}$.
\qed

As a consequence of this lemma we see that
Brian\c{c}on's Theorem is equivalent to saying that $X^{[n]}_{p,curv}$ is dense 
in $X^{[n]}_p$. This is a very important information: curvilinear subschemes 
are far easier to handle than any of the others. They contain only one subscheme
for any given smaller length, any small deformation of a curvilinear subscheme
is again locally curvilinear etc.

Generalising the definition of $X^{[n]}_p$ slightly, let $\Delta\subset S^nX$
denote the diagonal, and let $X^{[n]}_0:=\rho^{-1}(\Delta)$, endowed with
the reduced induced subscheme structure. Thus $X^{[n]}_0$ consists of all
subschemes $\xi\subset X$ of length $n$ which are supported at {\em some}
point in $X$. The fibres of the surjective
morphism $\rho:X_0^{[n]}\to X$ are the schemes $X^{[n]}_p$ considered above.
In fact, a choice of regular parameters near a point $p$ leads to a 
trivialisation of the morphism $\rho:X^{[n]}\to X$ near $p$, i.e.\ $\rho$
is a fibre bundle for the Zariski topology.

As an immediate consequence of Brian\c{c}on's Theorem we get

\begin{corollary}--- $X^{[n]}_0$ is an irreducible variety of dimension $n+1$.
\qed
\end{corollary}

Note that $X^{[n]}_p$ and $X^{[n]}_0$ have complementary dimensions as 
subvarieties in $X^{[n]}$. Their homological intersection is therefore zero-dimensional. However, the inclusion $X^{[n]}_p\subset X^{[n]}_0$ 
complicates the computation of the intersection product. The following
result was obtained by Ellingsrud and Str\o{}mme \cite{EllStrCoefficient}
by an inductive geometric argument:

\begin{theorem}[Ellingsrud, Str\o{}mme]\label{ESIntersectionnumber}---
$\deg([X^{[n]}_p]\cdot [X^{[n]}_0)=(-1)^{n-1}\cdot n$.\qed
\end{theorem}

 
\subsection{Incidence schemes}\label{subsectionIncidenceSchemes}

Since $X^{[n]}$ in facts represents the functor $\Hilb^n(X)$ of flat families
of subschemes of relative dimension 0 and length $n$, there is a {\em universal
family} of subschemes
$$\Xi_n\subset X^{[n]}\times X.$$
Again, for small values of $n$ there are explicit descriptions: $\Xi_0$ is empty, $\Xi_1$
is the diagonal in $X\times X$, and $\Xi_2$ is the blow-up
$\blowup{\Delta}{X\times X}$ of the diagonal in $X\times X$. The identification 
is given by the quotient map
$\blowup{\Delta}{X\times X}\to X^{[2]}=\blowup{\Delta}{X\times X}/\gothS_2$ 
and any of the two projections $\blowup{\Delta}{X\times X}\to X$.

Assume that $n'>n>0$. Then there is a uniquely determined closed 
subscheme $X^{[n',n]}\subset X^{[n']}\times X^{[n]}$ with the property that
any morphism $$f=(f_1,f_2):T\to X^{[n']}\times X^{[n]}$$
factors through
$X^{[n',n]}$ if and only if $f_{2,X}^{-1}(\Xi_n)\subset f_{1,X}^{-1}(\Xi_{n'})$. Closed points in $X^{[n',n]}$ correspond
to pairs $(\xi',\xi)$ of subschemes with $\xi\subset\xi'$. Let
$$X^{[n']}\stackrel{p_1}{\longleftarrow}X^{[n',n]}
\stackrel{p_2}{\longrightarrow}X^{[n]}$$
denote the two projections. Then $X^{[n',n]}$ parametrises two flat families
$$p_{2,X}^{-1}(\Xi_n)\subset p_{1,X}^{-1}(\Xi_{n'}).$$
Consider the corresponding exact sequence
\beeq{idealsheafInnsequence}
\ses{\ki_{n',n}}{p_{1,X}^*\ko_{\Xi_{n'}}}{p_{2,X}^*\ko_{\Xi_n}}.
\eneq

The ideal sheaf $\ki_{n',n}$ is a coherent sheaf on $X^{[n',n]}\times X$ which
is flat over $X^{[n',n]}$ and fibrewise zero-dimensional of length $n'-n$.
It therefore induces a classifying morphism to the symmetric product, analogously to the Hilbert-Chow morphism, which we will also denote by
$$\rho:X^{[n',n]}\to S^{n'-n}X.$$
As before let $X^{[n',n]}_0:=\rho^{-1}(\Delta)$, where $\Delta\subset S^{n'-n}X$
is the small diagonal. A point in $X^{[n',n]}_0$ is a triple
$(\xi',x,\xi)$ with $\xi\subset\xi'$ and $\Supp(\ki_{\xi/\xi'})=\{x\}$. 

We may decompose $X^{[n',n]}_0$ into locally closed subsets $Z_\ell$, $\ell\geq0$, with
$$Z_\ell:=\{(\xi',x,\xi)|\ell(\xi_x)=\ell\}.$$

\begin{lemma}\label{Z0andZ1}--- $Z_0$ and $Z_1$ are irreducible of dimension $n+n'+1$ and
$n+n'$, respectively, and $\dim(Z_\ell)<n+n'$ for all $\ell>1$. Moreover,
$Z_1$ is contained in the closure of $Z_0$.
\end{lemma}

\prf If $\ell=0$ or $1$,  the map $(\xi',x,\xi)\mapsto (\xi-\xi_x,\xi_x')$
is an open immersion 
$$Z_\ell\lra X^{[n-\ell]}\times X^{[n'-n+\ell]}_0.$$
It follows from Brian\c{c}on's Theorem that $Z_\ell$ is irreducible and 
$$\dim(Z_\ell)=2(n-\ell)+(n'-n+\ell+1)=n+n'+1-\ell.$$
For $\ell\geq 2$ consider the embedding
$$Z_\ell\lra X^{[n-\ell]}\times (X^{[\ell]}_0\times_XX^{[n'-n+\ell]}_0),
\quad(\xi',x,\xi)\mapsto(\xi-\xi_x, \xi_x,\xi'_x).$$
In fact, the image of $Z_\ell$ is contained in a {\em proper} closed subset of
the target variety: For {\em either} $\xi'_x$ is curvilinear, in which
case there is only a unique subscheme $\xi_x\subset \xi'_x$ of length $\ell$,
{\em or} $\xi'_x$ is not curvilinear and therefore contained in a proper
closed subset of $X^{[n'-n+\ell]}_0$. Now, the variety on the right hand side
has dimension $$2(n-\ell)+(\ell+1)+(n'-n+\ell+1)-2=n+n'.$$
 
Finally, a general point in $Z_1$ is of the form 
$(\zeta\cup\eta,x,\zeta\cup\{x\})$
where $\eta$ is a curvilinear subscheme supported at $x$ and 
disjoint from $\zeta$. Now it is easy to deform $\eta$ to a subscheme $\{x\}\cup\eta'$ with $\eta'$ supported at a point $x'\neq x$. Hence
a general point of $Z_1$ deforms into $Z_0$.\qed

\begin{definition}\label{definitionofQ}--- For any pair of nonnegative integers define subvarieties
$$E^{[n',n]},Q^{[n',n]}\subset X^{[n']}\times X\times X^{[n]}$$
as follows: if $n'>n>0$ let $Q^{[n',n]}$ and $E^{[n',n]}$ be the closure of
$Z_0$ and $Z_1$, respectively. Moreover, $Q^{[n',0]}:=X^{[n']}_0$, $E^{[n',0]}:=\emptyset$ and $Q^{[n,n]}:=\emptyset$, whereas $E^{[n,n]}:=\{
(\xi,x,\xi)|x\in\xi\}\isom \Xi_n$. On the other hand, if $n\geq n'$, let
$Q^{[n',n]}=T(Q^{[n,n']})$ and $E^{[n',n]}=T(E^{[n,n']})$ under the twist
$$T:X^{[n]}\times X\times X^{[n']}\to X^{[n']}\times X\times X^{[n]}.$$
\end{definition}

By construction $Q^{[n,n']}$ and $E^{[n,n']}$ are empty or irreducible varieties
of dimension $n+n'+1$ and $n+n'$, respectively.

Let us return to the particular case $n'-n=1$, the most basic of all incidence
situations: consider the projectivisation $\sigma:\IP(\ki_{\Xi_n})\to X^{[n]}\times X$. It is an easy exercise to see that there is a natural isomorphism $\IP(\ki_{\Xi_n})\isom X^{[n+1,n]}$ such that the diagram
$$\begin{array}{ccc}
\IP(\ki_{\Xi_n})&\rpfeil{5em}{\isom}&X^{[n+1,n]}\\
&{\scriptstyle\sigma}\searrow\quad\quad{\scriptstyle (p_2,\rho)}\swarrow\\[2ex]
&X^{[n]}\times X
\end{array}
$$
commutes.

The following theorem has independently been proved by Cheah \cite{CheahSmooth}, Ellings\-rud, and Tikhomirov (unpublished).

\begin{theorem}--- The incidence scheme $X^{[n+1,n]}$ is a smooth ir\-re\-duc\-ible variety.
\end{theorem}

An immediate corollary is the following: there is a natural closed immersion
$\blowup{\Xi_n}{X^{[n]}\times X}\to \IP(\ki_{\Xi_n})$; since both are 
irreducible varieties, this must be an isomorphism. The exceptional divisor 
$E$ is precisely the variety $E^{[n+1,n]}$ defined above. Hence in this situation we may write the sequence (\ref{idealsheafInnsequence}) as
\beeq{idealsheafgleichOminusE}
\ses{(\id,\rho)_*\ko_{X^{[n+1,n]}}(-E)}{p_{1,X}^*\ko_{\Xi_{n+1}}}{p_{2,X}^*
\ko_{\Xi_n}}.
\eneq
 

\section{The structure of the cohomology}\label{TheStructure}


The motivating problem in this study is to understand
the cohomology rings $H^*(X^{[n]})$ in terms of the cohomology ring $H^*(X)$.
For the symmetric product Grothendieck \cite{GrothSymmProd} showed that the 
natural map
$$\pi^*:H^*(S^nX;\IQ)\lra H^*(X^n;\IQ)\isom H^*(X;\IQ)^{\tensor n}$$
is an isomorphism onto the subring of invariant elements under the action of
$\gothS_n$. From this Macdonald \cite{Macdonald} computed the following formula for the 
Betti numbers of $S^nX$ by a purely algebraic argument:

\begin{theorem}[Macdonald]--- 
The Betti numbers of the symmetric products are given by the formula
$$\sum_{n\geq0}\sum_{i\geq0}b_i(S^nX)t^iq^n=
\prod_{i=0}^{2dim(Y)}(1-(-1)^it^iq)^{-(-1)^ib_i(X)}.$$\qed
\end{theorem}

For the Hilbert schemes the corresponding question for the Betti numbers is
much more difficult. This problem was solved by G\"ottsche \cite{Goettsche}:

\begin{theorem}[G\"ottsche]\label{GoettschesFormel}--- The Betti numbers $b_i(X^{[n]})$ are determined by the Betti numbers $b_j(X)$. More precisely,
the following formula holds:
$$\sum_{n\geq0}\sum_{i\geq0}b_i(X^{[n]})t^iq^n=\prod_{m>0}\prod_{j\geq0}(1-(-1)^{j}t^{2m-2+j}q^{m})^{-(-1)^jb_j(X)}$$
\end{theorem}

G\"ottsches original proof uses the Weil Conjectures \cite{Goettsche}. For a
different approach see \cite{Cheah}.

Among other things one learns from this formula that it is a good idea
to consider all Hilbert schemes simultaneously. This will become even more
striking through Nakajima's method which we will review in the next sections.
As a preparation we collect a few definitions:

\begin{definition}--- Let $\IH:=\bigoplus_{n,i\geq 0}\IH^{n,i}$ denote the
double graded vector space with components $\IH^{n,i}=H^i(X^{[n]};\IQ)$.
Since $X^{[0]}$ is a point, $\IH^{0,0}=\IQ$. The unit in $H^0(X^{[0]};\IQ)$ is called the `vacuum vector' and denoted by $\vacuum$.
\end{definition}

A linear map $\gothf:\IH\to\IH$ is homogeneous of bidegree $(\nu,\iota)$ if
$\gothf(\IH^{n,i})\subset\IH^{n+\nu,i+\iota}$ for all $n$ and $i$. If $\gothf,\gothf'\in \End(\IH)$
are homogeneous linear maps of bidegree $(\nu,\iota)$ and $(\nu',\iota')$, 
respectively, their commutator is defined by
$$[\gothf,\gothf']=\gothf\circ\gothf'-(-1)^{\iota\cdot\iota'}\gothf'\circ\gothf.$$
We use the notation $|\alpha|$, $|\gothf|$ etc.\ to denote
the cohomological degree of homogeneous cohomology classes, homogeneous linear
maps etc.

Setting
$$(\alpha,\beta):=\int_{X^{[n]}}\alpha\beta$$
for any $\alpha,\beta\in H^*(X^{[n]};\IQ)$ defines a non-degenerate 
(anti)symmetric bilinear form on $H^*(X^{[n]};\IQ)$ and hence on $\IH$.
For any homogeneous linear map $\gothf:\IH\to \IH$ its adjoint
$\gothf^\dagger$ is characterised by the relation
$$(\gothf(\alpha),\beta)=
(-1)^{|\gothf|\cdot|\alpha|}(\alpha,\gothf^\dagger(\beta)).$$
Clearly, $(\gothf\circ\gothg)^\dagger=\gothg^\dagger\circ\gothf^\dagger$.

\subsection{Correspondences}

Let $Y_1$ and $Y_2$ be smooth projective varieties, and let $u$ be a 
class in the Chow group $A_n(Y_1\times Y_2)$. (We tacitly assume rational
coefficients. This will not always be necessary. On the other hand, we are
not interested in integrality questions for the moment, and hence will not
pay attention to this problem). The image of $u$
in $H_{2n}(Y_1\times Y_2)$ will be denoted by the same symbol.
$u$ induces a homogeneous linear map
$$u_*:H^i(Y_2)\to H^{i+2(\dim Y_1-n)}(Y_1),\quad y\mapsto PD^{-1}p_{1*}(u\cap p_2^*y),
$$
where $PD:H^*(Y_1)\to H_*(Y_1)$ is the Poincar\'e duality map.

Assume that $Y_3$ is another smooth projective variety, and $v\in A_m(Y_2\times Y_3)$. Let $p_{ij}$ be the projection from $Y_1\times Y_2\times Y_3$ to
the  factors $Y_i\times Y_j$, and consider the element
$$w:=p_{13*}(p_{12}^*u\cdot p_{23}^*v)\in A_{n+m-\dim Y_2}(Y_1\times Y_3).
$$
Then
$$w_*=u_*\circ v_*.$$
See \cite[Ch.\ 16]{Fulton} for details. 

Suppose $U\subset Y_1\times Y_2$ and $V\subset Y_2\times Y_3$ are closed
subschemes such that $u\in A_*(U)$ and $v\in A_*(V)$. Let
$$W:=p_{13}(p_{12}^{-1}(U)\cap p_{23}^{-1}(V))$$
Then the class $w$ defined above is already defined in $A_*(W)$.

The following type of arguments will often show up in the sequel: one shows
that the dimension of $W$ is smaller than the degree of $w$, which forces
$w$ to be zero; or that there is at most one irreducible component $W_0$
of $W$ of maximal dimension with `correct' dimension $\dim(W_0)=\deg(w)$.
In this case one must have $w=\mu\cdot [W_0]$ and it suffices to determine
the multiplicity $\mu$.
 
Let $T:Y_1\times Y_2\to Y_2\times Y_1$ exchange the factors. Then a Chow cycle
$u$ induces two maps
$$u_*:H^*(Y_2)\to H^*(Y_1)\quad\mbox{and}\quad(Tu)_*:H^*(Y_1)\to H^*(Y_2)$$
which are related by the formula
$$\int_{Y_1}u_*(\alpha)\cdot\beta=\int_{Y_2}\alpha\cdot(Tu)_*(\beta).$$
This follows directly from the projection formula. Thus $(Tu)_*=u_*^\dagger$.

The following operators were introduced by Nakajima \cite{Nakajima}. The
study of their properties is the major theme of this article. We take the liberty to change the notations and sign conventions.

Recall that we defined (\ref{definitionofQ}) subvarieties
$$Q^{[n_1,n_2]}\subset X^{[n_1]}\times X\times X^{[n_2]}$$
of dimension $n_1+n_2+1$. Their fundamental classes are cycles
$$[Q^{[n_1,n_2]}]\in A_{n_1+n_2+1}(X^{[n_1]}\times X\times X^{[n_2]}).$$
Let the projections to the factors be denoted by $p_1$, $\rho$ and $p_2$.

\begin{definition}[Nakajima]\label{definitionofoscillators}--- Define
linear maps
$$\qq_{\ell}:H^*(X;\IQ)\lra \End(\IH),\qquad \ell\in\IZ,$$
as follows: assume first that $\ell\geq0$.
For $\alpha\in H^*(X;\IQ)$ and $y\in H^*(X^{[n]};\IQ)$ let
$$\qq_\ell(\alpha)(y):=[Q^{[n+\ell,n]}]_*(\alpha\tensor y)=PD^{-1}p_{1*}(
[Q^{[n+\ell,n]}]\cap (\rho^*\alpha\cdot p_2^*y)).$$
The operators for negative indices then are determined by the relation $$\qq_{-\ell}(\alpha):=(-1)^\ell\qq_{\ell}(\alpha)^\dagger.$$
\end{definition}

By definition, $\qq_\ell(\alpha)$ is a homogeneous linear map of
bidegree $(\ell, 2\ell-2+|\alpha|)$. Moreover, $\qq_0=0$, and if $\ell>0$,
the operator $\qq_{\ell}(\alpha)^\dagger$ is induced by the subvarieties
$Q^{[n,n+\ell]}$, $n\geq 0$.

\subsection{Nakajima's Main Theorem}

In this section we review the main result of \cite{Nakajima} and some of the
immediate consequences. Similar results have been announced by Grojnowski
\cite{Grojnowski}. 

\begin{theorem}[Nakajima]\label{oscillatorrelation}---
For any integers $n$ and $m$ and cohomology classes $\alpha$ and $\beta$,
the operators $\qq_n(\alpha)$  and $\qq_m(\beta)$ satisfy the following `oscillator relations':
$$[\qq_n(\alpha),\qq_m(\beta)]=n\cdot\delta_{n+m}\cdot\int_X\alpha\beta\cdot\id_{\IH}.$$\qed
\end{theorem}

Here and in the following we adopt the convention that $\delta_\nu$ equals
$1$ if $\nu=0$ and is zero else, and that any integral $\int_Z\alpha$ is zero
if $\deg(\alpha)\neq\dim_{\IR}(Z)$.

In \cite{Nakajima} Nakajima only showed that the commutator
relation holds with some universal nonzero constant instead of the coefficient
$n$. The correct value was computed directly by Ellingsrud and Str\o{}mme
\cite{EllStrCoefficient}: up to a sign factor, which depends on our
convention, this number is the intersection number of Theorem \ref{ESIntersectionnumber}. There is a different proof due to Grojnowski \cite{Grojnowski} and Nakajima \cite{NakajimaLectures} using `vertex operators'.

Consider the vector spaces
$$W_+:=H^*(X;\IQ)\tensor t\IQ[t]\quad\mbox{ and }\quad
W_-:=H^*(X;\IQ)\tensor t^{-1}\IQ[t^{-1}].$$ 
Define a non-degenerate skew-symmetric pairing on the vector space
$W:=W_-\oplus W_+$
by
$$\{\alpha\tensor t^n,\beta\tensor t^m\}
:=n\cdot\delta_{n+m}\cdot\int_X\alpha\beta.$$
Note that we are taking the expression `skew-symmetric' in a graded sense:   $$\{\alpha\tensor t^n,\beta\tensor t^m\}=-(-1)^{|\alpha|\cdot|\beta|}\{\beta\tensor t^m,\alpha\tensor t^n\}.$$
The {\em oscillator algebra} is the quotient of the tensor algebra $\kt W$ by
the two-sided ideal $I$ generated by the expressions $[v,w]-\{v,w\}\cdot 1$
with $v,w\in W$:
$$\kh:=\kt W/I.$$
$\kh$ is the (restricted) tensor product of countably many copies of Clifford
algebras arising from $H^{odd}(X;\IQ)$ and countably many copies of Weyl algebras arising from $H^{even}(X;\IQ)$.
As $W_+$ is isotropic with respect to the skew-form
$\{\,,\,\}$, the subalgebra in $\kh$ generated by $W_+$ is the symmetric algebra
$S^*W_+$ (taken again in a $\IZ/2$-graded sense). This becomes a double graded vector
space if we define the bidegree of $\alpha\tensor t^n$ as $(n,2n-2+|\alpha|)$.

Using these notations, Nakajima's Theorem can be rephrased by saying:\\
Sending $\alpha\tensor t^n\in W$ to $\qq_n(\alpha)\in\End(\IH)$ defines a
representation of $\kh$ on $\IH$.

The subspace $W_-$ of monomials of negative degree annihilates the vacuum
vector $\vacuum\in \IH$ for obvious degree reasons. Hence there is an embedding 
$$S^*W_+\isom \kh/\kh\cdot W_{-}\stackrel{\cdot\vacuum}{\lra} \kh\cdot\vacuum\subset \IH.$$
It is not difficult to check that the Poincar\'e series of $S^*W_+$ equals the
right hand side of G\"ottsche's formula. This implies:

\begin{corollary}[Nakajima]--- The action of $\kh$ on $\IH$ induces a module
isomorphism $S^*W_+\to \IH$. In particular, $\IH$ is irreducible and generated
by the vacuum vector.\qed
\end{corollary}


\section{The boundary operator}\label{TheDerivative}

The key to our solution of the Chern class problem is the introduction of
the boundary operator $\dd\in \End(\IH)$. This is done in \ref{TheBoundary}. We
begin with the discussion of related topics and ingredients for later proofs.

\subsection{Virasoro generators} 
 
Starting from the basic generators $\qq_n$ and the fundamental oscillator 
relations we will define the corresponding Virasoro generators $\LL_n$
in analogy to the procedure in conformal field theory. We will then give
concrete geometric interpretations for these generators.

Let $\delta:H^*(X)\to H^*(X\times X)=H^*(X)\tensor H^*(X)$ be the push-forward
map associated to the diagonal embedding. Equivalently, this is the linear map 
adjoint to the cup-product map.
If $\delta(\alpha)=\sum_i\alpha_i'\tensor\alpha_i''$, we will write
$\qq_n\qq_m\delta(\alpha)$ for $\sum_i\qq_n(\alpha_i')\qq_m(\alpha_i'').$

\begin{definition}--- Define operators $\LL_n:H^*(X;\IQ)\to\End(\IH)$,
$n\in\IZ$, as follows:  $$\LL_n:=\frac{1}{2}\sum_{\nu\in\IZ}\qq_\nu\qq_{n-\nu}\delta,\quad\mbox{ if } n\neq0$$
and
$$\LL_0:=\sum_{\nu>0}\qq_\nu\qq_{-\nu}\delta.$$
\end{definition}

\begin{remark}--- i)
The sums that appear in the definition are formally infinite.
However, as operators on any fixed vector in $\IH$, only finitely many of them
are nonzero. Hence the sums are locally finite and the operators $\LL_n$ 
are well-defined. $\LL_n(\alpha)$ is homogeneous of bidegree $(n,2n+|\alpha|)$
 
ii) Using the physicists' normal order convention
$$:\qq_n\qq_m:\,\,:=\left\{
\begin{array}{cl}
\qq_n\qq_m&\mbox{if }n\geq m,\\
\qq_m\qq_n&\mbox{if }n\leq m,
\end{array}\right.$$
the operators $\LL_n$ can be uniformly expressed as
$$\LL_n=\frac{1}{2}\sum_{\nu\in\IZ}:\qq_\nu\qq_{n-\nu}:\delta.$$
\end{remark}

\begin{theorem}\label{Lqcommutator}---
The operators $\LL_n$ and $\qq_m$ satisfy the following
commutation relations:
\begin{enumerate}
\item $[\LL_n(\alpha),\qq_m(\beta)]=-m\cdot\qq_{n+m}(\alpha\beta).$
\item $[\LL_n(\alpha),\LL_m(\beta)]=(n-m)\cdot\LL_{n+m}(\alpha\beta)-\frac{n^3-n}{12}\delta_{n+m}\cdot\int_Xc_2(X)\alpha\beta\cdot\id_{\IH}.$
\end{enumerate}
\end{theorem}

Taking only the operators $\LL_n(1)$, $n\in IZ$, we see that the Virasoro 
algebra acts on $\IH$ with central charge equal to the Euler number of $X$.

\prf Assume first that $n\neq 0$. For any classes $\alpha$ and $\beta$ with 
$$\delta(\alpha)=\sum_i\alpha'_i\tensor\alpha''_i$$
we have
\begin{eqnarray*}
[\qq_\nu(\alpha'_i)\qq_{n-\nu}(\alpha''_i),\qq_m(\beta)]&=&
\qq_\nu(\alpha'_i)[\qq_{n-\nu}(\alpha''_i),\qq_m(\beta)]\\
&&+(-1)^{|\beta|\cdot|\alpha''_i|}[\qq_\nu(\alpha'_i),\qq_m(\beta)]
\qq_{n-\nu}(\alpha''_i)\\
&=&(-m)\delta_{n+m-\nu}\cdot\qq_{n+m}(\alpha'_i)\cdot\int_X\alpha''_i\beta\\
&&+(-1)^{|\beta|\cdot|\alpha|}(-m)\delta_{\nu+m}\cdot\int_X\beta\alpha'_i\cdot\qq_{n+m}(\alpha_i'').
\end{eqnarray*}
If we sum up over all $\nu$ and $i$, we get
$$2[\LL_n(\alpha),\qq_m(\beta)]=\sum_{\nu}
[\qq_\nu\qq_{n-\nu}\delta(\alpha),\qq_m(\beta)]=(-m)\cdot\qq_{n+m}(\gamma)
$$
with
$$\gamma=pr_{1*}(\delta(\alpha)\cdot pr_2^*(\beta))+(-1)^{|\beta|\cdot|\alpha|}
\cdot pr_{2*}(pr_1^*(\beta)\cdot\delta(\alpha))=2\cdot \alpha\beta.$$
Similarly, for $\nu>0$,
$$
[\qq_{\nu}\qq_{-\nu}\delta(\alpha),\qq_m(\beta)]=-m\cdot\qq_m(\alpha\beta)\cdot
(\delta_{m-\nu}+\delta_{m+\nu}).$$
Thus summing up over all $\nu>0$ we find again 
$$[\LL_0(\alpha),\qq_m(\beta)]=-m\cdot \qq_m(\alpha\beta).$$
This proves the first part of the theorem. 

As for the second part, assume first that $n\geq0$. In order to avoid
case considerations let us agree that $\qq_{\frac{N}{2}}$ is zero if $N$ is odd. Then
we may write:
$$\LL_m=\frac{1}{2}\qq_{\frac{m}{2}}^2\delta+\sum_{\mu>\frac{m}{2}}\qq_\mu\qq_{m-\mu}\delta.$$
By the first part of the theorem we have
$$[\LL_n(\alpha),\qq_\mu\qq_{m-\mu}\delta(\beta)]
=\Big(-\mu\qq_{n+\mu}\qq_{m-\mu}+(\mu-m)\qq_\mu\qq_{n+m-\mu}\Big)
\delta(\alpha\beta).$$
In the following calculation we suppress $\alpha,\beta$ and $\delta$ up to the very
end. Summing up over all $\mu\geq 0$, we get:
\begin{eqnarray*}
[\LL_n,\LL_m]&=&-\frac{m}{4}(\qq_{n+\frac{m}{2}}\qq_{\frac{m}{2}}+\qq_{\frac{m}{2}}\qq_{n+\frac{m}{2}})\\
&&+\sum_{\mu>\frac{m}{2}}(\mu-m)\qq_{\mu}\qq_{n+m-\mu}
+\sum_{\mu>\frac{m}{2}}(-\mu)\qq_{n+\mu}\qq_{m-\mu}\\
&=&-\frac{m}{4}(\qq_{n+\frac{m}{2}}\qq_{\frac{m}{2}}+\qq_{\frac{m}{2}}\qq_{n+\frac{m}{2}})\\
&&+\sum_{\mu>\frac{m}{2}}(\mu-m)\qq_{\mu}\qq_{n+m-\mu}
+\sum_{\mu>n+\frac{m}{2}}(n-\mu)\qq_{\mu}\qq_{n+m-\mu}
\end{eqnarray*}
Hence
\begin{eqnarray*}
[\LL_n,\LL_m]-(n-m)\sum_{\mu>\frac{n+m}{2}}\qq_\mu\qq_{n+m-\mu}&=&
-\frac{m}{4}(\qq_{n+\frac{m}{2}}\qq_{\frac{m}{2}}+\qq_{\frac{m}{2}}\qq_{n+\frac{m}{2}})\\
&&+\sum_{\frac{m}{2}<\mu\leq\frac{m+n}{2}}(\mu-m)\qq_\mu\qq_{m+n-\mu}\\
&&-\sum_{\frac{n+m}{2}<\mu\leq n+\frac{m}{2}}(n-\mu)\qq_\mu\qq_{n+m-\mu}
\end{eqnarray*}
Now split off the summands corresponding to the indices $\mu=\frac{m+n}{2}$
and $\mu=n+\frac{m}{2}$ from the sums. Substituting $n+m-\mu$
for $\mu$ in the second sum on the right hand side, we are left with the
expression:
$$[\LL_n,\LL_m]-(n-m)\LL_{n+m}=
-\frac{m}{4}[\qq_{\frac{m}{2}},\qq_{n+\frac{m}{2}}]+\sum_{\frac{m}{2}<\mu<\frac{n+m}{2}}(\mu-m)[\qq_\mu,\qq_{n+m-\mu}]
$$
The right hand side is zero unless $n+m=0$. In this case, observe that the
composition
$$H^*(X)\stackrel{\delta}{\lra}H^*(X)\tensor H^*(X) \stackrel{\cup}\lra H^*(X)$$
is multiplication with $c_2(X)$. Hence we see that
$$[\LL_n(\alpha),\LL_m(\beta)]=(n-m)\LL_{n+m}(\alpha\beta)+\delta_{n+m}\cdot\int_Xc_2(X)\alpha\beta\cdot N,$$
where $N$ is the number
$$N=\sum_{0<\nu<\frac{n}{2}}\nu(\nu-n)\qquad\mbox{ if }n\mbox{ is odd,}$$
and 
$$N=\sum_{0<\nu<\frac{n}{2}}\nu(\nu-n)-\frac{n^2}{8}\quad\mbox{ if }n\mbox{ is even.}
$$
An easy computation shows that in both cases $N$ equals $(n-n^3)/12$.
\qed

Recall the definition of the varieties
$E^{[n,n']}\subset X^{[n]}\times X\times X^{[n']}$ in (\ref{definitionofQ}).

\begin{definition}--- Let $\ell$ be a nonnegative integer and let $$\ee_\ell:H^*(X)\to \End(\IH)$$
be the linear map
$$\ee_\ell(\alpha)(y)=[E^{[n+\ell,n]}]_*(\alpha\tensor y)=PD^{-1}p_{1*}(
[E^{[n+\ell,n]}]\cap(\rho^*\alpha\cdot p_{2}^*y))$$
for $\alpha\in H^*(X;\IQ)$ and $y\in H^*(X^{[n]};\IQ)$.
\end{definition}

The following theorem gives a `finite' geometric interpretation of the infinite
sums which define the Virasoro operators.

\begin{theorem}\label{eqcommutator}--- Let $n$ be a nonnegative integer.
\begin{enumerate}
\item
$$[\ee_n(\alpha),\qq_m(\beta)]=\left\{
\begin{array}{ll}
m\cdot\qq_{n+m}(\alpha\beta)&\mbox{if }m>0\mbox{ or }m<-n.\\
0&\mbox{else}.
\end{array}
\right.$$
\item
$$\ee_n+\LL_n=\frac{1}{2}\sum_{0<\nu<n}\qq_\nu\qq_{n-\nu}\delta.$$ 
\end{enumerate}
\end{theorem}

\prf Ad 1: Assume first that $m\geq1$. To simplify the notations we introduce
the short-hand
$$X^{[n_1],[n_2],\ldots,[n_k]}:=X^{[n_1]}\times X^{[n_2]}\times \ldots\times X^{[n_k]}$$
Suppose $\ell\geq0$, and consider the following diagram
$$\begin{array}{ccccc}
X^{[\ell+n+m],[1],[\ell+m]}&\stackrel{p_{123}}{\verylongleftarrow{3em}}&
X^{[\ell+n+m],[1],[\ell+m],[1],[\ell]}&
\stackrel{p_{345}}{\verylongrightarrow{3em}} X^{[\ell+m],[1],[\ell]} \\
&&\phantom{\scriptstyle p_{1245}}\Bigg\downarrow
{\scriptstyle p_{1245}}\\[3ex]
&&X^{[\ell+n+m],[1],[1],[\ell]}
\end{array}$$
The product operator $\ee_n\qq_m$ is induced by the class
$$z:=p_{1245*}(p_{123}^*[E^{[\ell+m+n,\ell+m]}]\cdot p_{345}^*[Q^{[\ell+m,\ell]}])\in A_{2\ell+n+m+1}(Z')$$
where
\begin{eqnarray*}Z'&:=&p_{1245}(p_{123}^{-1}(E^{[\ell+m+n,\ell+m]})\cap
p_{345}^{-1}(Q^{[\ell+m,\ell]}))\\
&\subset&Z:=\{(\xi',x,y,\xi)|\exists\eta:\xi'-\eta=nx,\eta-\xi=my, x\in\eta\}
\end{eqnarray*}
Here the notation $\eta-\xi=my$ should comprise the conditions: $\xi$ is 
a subscheme of $\eta$, and the ideal sheaf of $\xi$ in $\eta$ is of length $m$
and is supported at $y$ etc.

Similarly, the operator $\qq_n\ee_m$ is induced by a class
$v\in A_{2\ell+m+n+1}(V')$ with
$$V'\subset V:=\{(\xi',x,y,\xi)|\exists \eta':\xi'-\eta'=mx, \eta'-\xi=ny, y\in\xi\}.$$
Moreover, if $T:X^{[\ell+m+n],[1],[1],[\ell]}\lra X^{[\ell+m+n],[1],[1],[\ell]}$
exchanges the two copies of $X$ in the middle, then the commutator
$[\ee_n,\qq_m]$ is induced by $z-T(v)$.

Now observe that off the diagonal $\{x=y\}\subset X^{[\ell+m+n],[1],[1],[\ell]}$
the subsets $Z$ and $T(V)$ are equal. Moreover, there is only one component
of (maximal possible) dimension $2\ell+n+m+1$. It is easy to see that this
component has multiplicity $1$ both in $z$ and $T(v)$: the intersection
$$p_{123}^{-1}(E^{[\ell+m+n,\ell+m]})\cap
p_{345}^{-1}(Q^{[\ell+m,\ell]})$$
is transversal over a general point in this component of $Z$, and maps injectively into $Z$. 
Thus the only contributions to $z-T(v)$ may arise from 
the diagonal part. Now 
$$V\cap\{x=y\}=\{(\xi',x,x,\xi)|\xi'-\xi=(n+m)x, x\in \xi\}.$$
We have seen earlier (\ref{Z0andZ1})
that this set has dimension $\leq 2\ell+n+m$ and hence may be disregarded.
On the other hand
$$Z\cap\{x=y\}=\{(\xi',x,x,\xi)|\xi'-\xi=(n+m)x\}.$$
Again using \ref{Z0andZ1} we see that this set has only one component $D$ of
(maximal) dimension $2\ell+n+m+1$. Moreover, this component is the
image of the embedding
$$\iota:Q^{[\ell+n+m,\ell]}\to X^{[\ell+n+m],[1],[1],[\ell]}, (\xi',x,\xi)\mapsto (\xi',x,x,\xi).$$

Let $\alpha,\beta\in H^*(X;\IQ)$ and $y\in H^*(X^{[\ell]};\IQ)$. Then we have
\begin{eqnarray*}
\lefteqn{p_{1*}([D]\cap p_{23}^*(\alpha\tensor\beta)\cdot
p_4^* y)
}\hspace{3em}\\
&=&
p_{1*}(\iota_*[Q^{[\ell+n+m,\ell]}]\cap p_{23}^*(\alpha\tensor\beta)\cdot
p_4^* y)\\
&=&
p_{1*}([Q^{[\ell+n+m,\ell]}]\cap \iota^*(p_{23}^*(\alpha\tensor\beta)\cdot p_4^*y))\\
&=&p_{1*}([Q^{[\ell+n+m,\ell]}]\cap p_2^*(\alpha\beta)\cdot p_3^*y)
\end{eqnarray*}
This shows that 
$$[\ee_n(\alpha),\qq_m(\beta)]=\mu\cdot\qq_{n+m}(\alpha\beta)$$
for some integer $\mu$.
Hence it remains to compute the multiplicity $\mu$ of $[D]$ in $z$.
To this end we pick a general point $d\in D$ and inspect the 
intersection of $p_{123}^{-1}(E^{[\ell+n+m,\ell]})$ and
$p_{345}^{-1}(Q^{[\ell+m,\ell]})$ along the fibre $p_{1245}^{-1}(d)$.

A general point in $D$ is of the form
$$d=(\xi',x,x,\xi)\quad\mbox{ with }\quad\xi'=\xi\cup\zeta,$$
where $\zeta$ is a curvilinear subscheme of $X$ of length $n+m$, supported
in a single point $x$ which is disjoint from $\xi$.
Since $\zeta$ is curvilinear, there is a unique subscheme $\eta\subset\zeta$
of length $m$, and hence $p_{1245}^{-1}(d)$ consists of the single point
$$d'=(\xi\cup\zeta,x,\xi\cup\eta,x,\xi)$$
Near $d'$ the varieties $X^{[\ell+m+n],[1],[\ell+m],[1],[\ell]}$ and 
$X^{[\ell],[\ell],[\ell]}\times X^{[m+n],[1],[m],[1]}$ are locally isomorphic
in the \'etale topology;
and similarly $E^{[\ell+m+n,\ell+m]}$ to $X^{[\ell]}\times E^{[m+n,m]}$ and
$Q^{[\ell+m,\ell]}$ to $X^{[\ell]}\times X^{[m]}_0$. Thus we may split off the
factors $X^{[\ell]}$ from the geometric picture. In the end this amounts
to saying that we may assume without loss of generality that $\ell=0$.

Moreover, the calculation is local (in the \'etale topology) in $X$, so that we may assume that $X=\IA^{2}
=\Spec\IC[z,w]$ and $\ki_\zeta=(w,z^{n+m})$, $\ki_{\eta}=(w,z^m)$ and $\ki_x=(w,z)$. Then $d'$ has an affine neighbourhood $\isom\IA^{4m+2n+4}$
in $X^{[n+m],[1],[m],[1]}$ with coordinate functions
$$a_0,\ldots,a_{n+m-1},b_0,\ldots,b_{n+m-1},w_1,z_1,c_0,\ldots,c_{m-1},d_0,
\ldots,d_{m-1},w_2,z_2,$$
which parametrises quadruples $(\zeta,x,\eta,y)$ of subschemes in $X$ given by the ideals
$$(w-g_1(z),f_1(z)),\quad (w-w_1,z-z_1),\quad (w-g_2(z),f_2(z)),\quad (w-w_2,z-z_2),$$
where
$$f_1(z)=\sum_{i=0}^{n+m-1}a_iz^i+z^{n+m},\quad
g_1(z)=\sum_{i=0}^{n+m-1}b_iz^i
$$
and
$$f_2(z)=\sum_{i=0}^{m-1}c_iz^i+z^{m},\quad
g_2(z)=\sum_{i=0}^{m-1}d_iz^i.$$
Now $(\eta,y)$ belongs to $X^{[m]}_0$, i.e.\ $\Supp(\eta)=\{y\}$, if and only if
\beeq{conditionA}
f_2(z)=(z-z_2)^m\, \mbox{ and }\, w_2=g_2(z_2).
\eneq
And $(\zeta,x,\eta)$ belongs to $Q^{[n+m,m]}$ if and only if the following
three conditions are satisfied: $\eta\subset\zeta$, i.e.
\beeq{conditionB}
g_1(z)=g_2(z)+f_2(z)\cdot h(z)\,\mbox{ and }\,f_1(z)=f_2(z)\cdot k(z)
\eneq
with polynomials $h$ and $k$ of degree $n-1$ and $n$, respectively;
the ideal sheaf $\ki_{\eta/\zeta}$ is supported at $x$, i.e.\
\beeq{conditionC}
k(z)=(z-z_1)^m\,\mbox{ and }\, w_1=g_1(z_1)
\eneq
and finally, $x$ must be contained in $\eta$, which imposes the condition
\beeq{conditionD}
f_2(z_1)=0
\eneq
One easily checks that the equations (\ref{conditionA}) - (\ref{conditionC})
cut out a smooth subvariety which projects isomorphically to the affine space
$\Spec\,\IC[z_1,z_2,b_0,\ldots,b_{n+m-1}]$. Moreover,
in these coordinates the last condition (\ref{conditionD}) simply reads  $(z_1-z_2)^m=0$. Hence the multiplicity $\mu$ equals the exponent $m$.   

Next, we consider the case $[\ee_n,\qq_{-m}]$ with $0\leq m\leq n$. There is
nothing to prove if $m=0$. Hence assume that $m>0$. Dimension arguments similar
to the ones above show that the cycle $v$ which induces the commutator 
$[\qq_{-m},\ee_n]$ must be supported on the closed subsets
$$V:=\{(\xi,x,x,\zeta)|\xi\supset\zeta\ni x, \xi-\zeta=(n+m)x\}\subset
X^{[\ell+n-m],[1],[1],[\ell]},\quad\ell\geq0.$$
The cycle $v$ has degree $2\ell+n-m+1$, so that it suffices to show that
$\dim(V)\leq 2\ell+n-m$. This follows from Lemma \ref{Z0andZ1}.

It remains to consider the case $[\ee_n,\qq_m]$ with $m<-n$. A dimension check
of the set-theoretic support of the intersection cycle shows that we must
have
$$[\ee_n(\alpha),\qq_{-m}(\beta)]=\mu\cdot\qq_{n-m}(\alpha\beta)$$
for some integer $\mu$, independently of $\alpha$ and $\beta$. 
To determine $\mu$, we proceed algebraically and take the commutator with $\qq_{m-n}(1)$:
$$
[\,[\ee_n(\alpha),\qq_{-m}(\beta)],\qq_{m-n}(1)]=
\mu\cdot[\qq_{n-m}(\alpha\beta),\qq_{m-n}(1)]=
\mu(n-m)\int_X\alpha\beta\cdot\id_\IH.
$$
On the other hand, combining the Jacobi identity, the oscillator relations and 
the first part of the proof yields
\begin{eqnarray*}
[\,[\ee_n(\alpha),\qq_{-m}(\beta)],\qq_{m-n}(1)]&=&
[\,[\ee_n(\alpha),\qq_{m-n}(1)],\qq_{-m}(\alpha)]\\
&=&(m-n)[\qq_m(\alpha),\qq_{-m}(\beta)]\\
&=&m(m-n)\int_X\alpha\beta\cdot\id_{\IH}.
\end{eqnarray*}
It follows that $\mu=-m$.
\medskip

Ad 2:
Consider the difference 
$\gothy:=\ee_n(\alpha)+\LL_n(\alpha)-\frac{1}{2}\sum_{\nu=1}^{n-1}\qq_\nu\qq_{n-\nu}\delta(\alpha)$.
Comparing the expressions in \ref{Lqcommutator} and part 1 of the theorem
we see that $\gothy$ commutes with all operators $\qq_m$ , $m\in\IZ$. 
Since $\IH$ is a simple $\kn$-module, $\gothy$ must be a scalar (in some
algebraic extension of $\IQ$), which is impossible: if $n>0$, then $\gothy$ has non-trivial bidegree $(n,2n+|\alpha|)$, and if $n=0$, it is easy to see 
directly that $\gothy\cdot\vacuum=0$.
\qed
 
\begin{remark}--- In particular, the operator $\LL_0(\alpha)$ has the
following geometric interpretation: the universal family $\Xi_n\subset
X^{[n]}\times X$ induces a homomorphism
$$[\Xi_n]_*:H^*(X;\IQ)\lra H^*(X^{[n]};\IQ),$$
and
$$\LL_0(\alpha)(y)=-[\Xi_n]_*(\alpha)\cdot y \quad\mbox{ for all }\quad
y\in H^*(X^{[n]};\IQ).$$
If we insert $\alpha=-1_X$, we get $\LL_0(-1_X)(y)=n\cdot y$ for all $y\in H^*(X^{[n]};\IQ)$.
Thus $\LL_0(-1_X)$ is the `number' operator, that counts with how many points
we are dealing. This can, of course, also be deduced directly from the 
definition of $\LL_0$.
\end{remark}

\subsection{The boundary of the Hilbert scheme}\label{TheBoundary}

For any partition $\lambda=(\lambda_1\geq\lambda_2\geq\ldots\geq\lambda_s>0)$
of $n$ the tuples $\sum_{1\leq i\leq s}\lambda_ix_i$, $x_i\in X$, form
a locally closed subset $S^n_\lambda X$ in $S^nX$. Let
$X^{[n]}_\lambda=\rho^{-1}(S^n_\lambda X)$. It follows from Brian\c{c}on's
Theorem that $X^{[n]}_\lambda$ is irreducible and
$$\dim(X^{[n]}_{\lambda})=\sum_{1\leq i\leq s}(\lambda_i+1)=n+s.$$
The generic open stratum is 
$X^{[n]}_{(1,1,\ldots,1)}$. It corresponds to the configuration space of 
unordered $n$-tuples of pairwise distinct points. Furthermore, there is 
precisely one stratum of codimension 1, namely $X^{[n]}_{(2,1,\ldots,1)}.$

If $\lambda=(\lambda_1,\ldots,\lambda_s)$ and $\mu=(\mu_1,\ldots,\mu_{s'})$
are partitions of $n$, then $X^{[n]}_\mu$ is contained in the closure of $X^{[n]}_\lambda$ if and only if there is a surjection
$$\phi:\{1,\ldots,s\}\to \{1,\ldots,s'\}$$
such that $\mu_j=\sum_{i\in\phi^{-1}(j)}\lambda_i$ for all $j$.
It follows that
$$\partial X^{[n]}:=\bigcup_{\lambda\neq(1,\ldots,1)}X^{[n]}_\lambda=
\overline{X^{[n]}_{(2,1,\ldots,1)}}$$
is an irreducible divisor in $X^{[n]}$. As it is the complement of the configuration space in $X^{[n]}$ we might and will call it the {\em boundary}
of $X^{[n]}$.

We will need a different description of the divisor $\partial X^{[n]}$
in sheaf theoretic terms.
Let $p:\Xi_n\to X^{[n]}$ be the projection, and define sheaves
$$\ko_X^{[n]}:=p_*(\ko_{\Xi_n})\in \Coh(X^{[n]}).$$
As $p$ is flat and finite of degree $n$, $\ko_X^{[n]}$ is locally free of rank
$n$.  

\begin{lemma}\label{DnX}\label{boundarydescription}--- We have
$$\left[\partial X^{[n]}\right]=-2\,c_1(\ko_X^{[n]})$$ 
Moreover, let $E\subset X^{[n+1,n]}$ be the exceptional divisor. Then
$$p_1^*\partial X^{[n+1]}-p_2^*\partial X^{[n]}=2\cdot E.$$
\end{lemma}

\prf $\partial X^{[n]}$ is the branching divisor of the finite flat
morphism $\Xi_n\to X^{[n]}$. The assertion holds true in a more general setting: if $Y$ is a smooth variety and
$\pi:Y'\to Y$ is a finite flat map, so that $\ka:=\pi_*\ko_{Y'}$ is a locally
free $\ko_Y$--sheaf, the branching divisor is given by the discriminant of
the $\ko_Y$-bilinear form
$$\ka\tensor_{\ko_Y}\ka\stackrel{\cdot}\lra\ka\stackrel{\mbox{tr}}\lra\ko_Y,$$
or, equivalently, by the determinant of the adjoint linear map $\ka\to 
\ka\dual$, so that indeed the branching divisor is given by $-2\,c_1(\ka)$. 

Applying $p_*$ to the short exact sequence (\ref{idealsheafgleichOminusE}) we
get an exact sequence
$$\ses{\ko_{X^{[n+1,n]}}(-E)}{p_1^*\ko_X^{[n+1]}}{p_2^*\ko_X^{[n]}},$$
from which one deduces the second assertion.\qed

This proof was communicated to me by S.\ A.\ Str\o{}mme and replaces a slightly
longer one in an earlier version.

\begin{definition}---  Let $\dd:\IH\to \IH$ be the homogeneous linear map
of bidegree $(0,2)$ given by
$$\dd(x):=c_1(\ko_X^{[n]})\cdot x=-\frac{1}{2}\left[
\partial X^{[n]}\right]\cdot x\quad\mbox{ for all } x\in H^*(X^{[n]}).$$
For any endomorphism $\gothf\in\End(\IH)$ its derivative is  
$\gothf':=[\dd,\gothf]$. As usual, we write $\gothf^{(n)}:=(\ad\dd)^n(\gothf)$
for the higher derivatives.
\end{definition}

It follows directly from the Jacobi identity
that $\gothf\mapsto \gothf'$ is a derivation, i.e.\ for any two operators $\gotha,\gothb\in \End(\IH)$ the `Leibniz rule' holds:
$$(\gotha\gothb)'=\gotha'\gothb+\gotha\gothb'\quad\mbox{ and }\quad 
[\gotha,\gothb]'=[\gotha',\gothb]+[\gotha,\gothb'].$$
Moreover, if $\gothf:H^*(X^{[\ell]})\to H^*(X^{[n]})$ is a homogeneous linear
map, then $|\gothf'|=|\gothf|+2$, so that $\gothf$ and $\gothf'$ have the same
parity. Furthermore, 
$$(\gothf')^\dagger=-(\gothf^\dagger)'.$$
Indeed, this follows formally from the obvious fact that $\dd^\dagger=\dd$. 

Let $n'>n$ be nonnegative integers, and consider the incidence variety
$X^{[n',n]}\subset X^{[n']}\times X^{[n]}$. Recall the definition of the
ideal sheaf $\ki_{n',n}$ and the exact sequence
$$\ses{\ki_{n',n}}{p_{1,X}^*\ko_{\Xi_{n'}}}{p_{2,X}^*\ko_{\Xi_n}}.$$
Then $p_*(\ki_{n',n})$ is a locally free sheaf of rank $n'-n$ on $X^{[n',n]}$.

\begin{lemma}\label{wasistdieableitung}--- Let $u_*:H^*(X^{[n]};\IQ)\to 
H^*(X^{[n']};\IQ)$ be the induced linear map associated to a class
$u\in A_*(X^{[n',n]})$. Then
$$(u_*)'=(c_1(p_*(\ki_{n',n}))\cdot u)_*.$$
\end{lemma}

\prf Let $y\in H^*(X^{[n]};\IQ)$. Then
\begin{eqnarray*}
(u_*)'(y)&=&\dd(u_*(y))-u_*(\dd(y))\\
&=&
c_1(p_*\ko_{\Xi_{n'}})\cdot PD^{-1}p_{1*}(u\cdot p_2^*y)\\
&&-PD^{-1}p_{1*}(u\cdot p_2^*(c_1(p_*\ko_{\Xi_n})\cdot y))\\
&=&PD^{-1}p_{1*}((p_1^*c_1(p_*\ko_{\Xi_{n'}})-p_2^*c_1(p_*\ko_{\Xi_n}))\cdot u\cdot
p_2^*y)\\
&=&v_*(y)
\end{eqnarray*}
with $v=(p_1^*c_1(p_*\ko_{\Xi_{n'}})-p_2^*c_1(p_*\ko_{\Xi_n}))\cdot u$,
and
\begin{eqnarray*}
p_1^*c_1(p_*\ko_{\Xi_{n'}})-p_2^*c_1(p_*\ko_{\Xi_n})&=&c_1(p_*p_{1,X}^*\ko_{\Xi_{x'}})-c_1(p_*p_{2,X}^*\ko_{X_n})\\
&=&c_1(p_*\ki_{n',n}).
\end{eqnarray*}
\qed

\subsection{The derivative of $\qq_n$}

In order to understand the intersection behaviour of the boundary $\partial
X^{[n]}$ we need to know how the operator $\dd$ commutes with the 
basic operators $\qq_n$, in other words: we need to compute the derivative of 
$\qq_n$.

The following theorem describes the derivative of the operator $\qq_n$ in two ways:
By its action on any of the other basic operators, and as a polynomial
expression in the basic operators.

Let $K$ denote the canonical class of the surface $X$.

\begin{theorem}\label{maintheorem}--- For all  $n,m\in\IZ$ and  $\alpha,\beta\in H^*(X;\IQ)$ the following holds:
\begin{enumerate}
\item
$[\qq_n'(\alpha),\qq_m(\beta)]=-nm\cdot\left\{\qq_{n+m}(\alpha\beta)
+\frac{|n|-1}{2}\delta_{n+m}\cdot\int_XK\alpha\beta\cdot\id_{\IH}\right\}.$
\item
$\qq'_n(\alpha)=n\cdot\LL_n(\alpha)+\frac{n(|n|-1)}{2}\qq_n(K\alpha).$
\end{enumerate}
\end{theorem}

\begin{corollary}\label{q1machtsallein}--- The operators $\dd$ and $\qq_1(\alpha)$, $\alpha\in H^*(X)$, suffice to generate $\IH$ from the 
vacuum $\vacuum$.\qed
\end{corollary}

\prfofthetheorem The second assertion is an immediate consequence of the first:
by Nakajima's relations \ref{oscillatorrelation} and the relations
\ref{Lqcommutator} we see that
\begin{eqnarray*}\lefteqn{[n\cdot\LL_n(\alpha)+\frac{n(|n|-1)}{2}\qq_n(K\alpha),\qq_m(\beta)]=}
\hspace{5em}\\
&&-nm\cdot\qq_{n+m}(\alpha\beta)+\delta_{n+m}\frac{n^2(|n|-1)}{2}\int_{X}K\alpha\beta\cdot\id_{\IH}.
\end{eqnarray*}
Hence the difference of $\qq'_n$ and the expression on the right hand side
in the theorem commutes with all operators $\qq_m$, $m\in \IZ$. Since $\IH$ is 
an irreducible $\kn$-module, it follows from Schur's Lemma that this difference
is given by multiplication with a scalar
(say, after passage to some algebraic closure of $\IQ$). But this
is impossible for degree reasons: the bidegree of $\qq_n'(\alpha)$ is
$(n,2n+|\alpha|)$. (The case $n=0$ being trivial anyhow.)

The proof of the first assertion has two parts of quite different nature:
We need to distinguish the cases $n+m\neq0$ and $n+m=0$ and deal with them
separately.

\begin{proposition}\label{easycase}--- $[\qq_n'(\alpha),\qq_m(\beta)]=-nm\cdot\qq_{n+m}(\alpha\beta)$ for any
two integers $n,m$ with $n+m\neq0$ and cohomology classes $\alpha,\beta
\in H^*(X)$.
\end{proposition}

\prf {\em Step 1:} Assume that $n$ and $m$ are positive.
We proceed as in the proof of Theorem \ref{eqcommutator}. Let $\ell$ be
nonnegative, and consider the diagram
$$\begin{array}{ccccc}
X^{[\ell+n+m],[1],[\ell+m]}&\stackrel{p_{123}}{\verylongleftarrow{3em}}&
X^{[\ell+n+m],[1],[\ell+m],[1],[\ell]}&
\stackrel{p_{345}}{\verylongrightarrow{3em}} X^{[\ell+m],[1],[\ell]} \\
&&\phantom{\scriptstyle p_{1245}}\Bigg\downarrow
{\scriptstyle p_{1245}}\\[3ex]
&&X^{[\ell+n+m],[1],[1],[\ell]}.
\end{array}$$
Let 
$$v:=p_{123}^*[Q^{[\ell+m+n,\ell+m]}]\cdot p_{345}^*[Q^{[\ell+m,\ell]}]\in
A_{2\ell+m+n+2}(V),$$
{}$$V:=p_{123}^{-1}(Q^{[\ell+m+n,\ell+m]})\cap 
p_{345}^{-1}(Q^{[\ell+m,\ell]}).$$
According to Lemma \ref{wasistdieableitung}, the operator $\qq_n'\qq_m$
is induced by the class
$$w=p_{1245*}(p_{123}^*c_1(\ki_{\ell+m+n,\ell+m})\cdot v)\in A_{2\ell+m+n+1}(W), W:=p_{1245}(V).$$
Let $V'\subset V$ and $W'\subset W$ denote the open subsets of those tuples
$(\xi,x,\sigma,y,\zeta)$ and $(\xi,x,y,\zeta)$, respectively, where either
$x\neq y$ or $x=y$ but $\xi_x$ is curvilinear. Certainly, $V'=p_{1245}^{-1}(W')$, but in fact we even have that 
$p_{1245}:V'\to W'$ is an isomorphism: for the conditions imposed on $V'$ imply
that $\sigma$ is already determined by the remaining data $(\xi,x,y,\zeta)$.

{\em Claim: $V'$ is irreducible of dimension $2\ell+n+m+2$.}

For it follows from Brian\c{c}on's Theorem that the open part $V'\setminus\{x=y\}$ is irreducible of dimension $2\ell+(n+1)+(m+1)$, and
tuples of the second kind, i.e.\ $(\xi,x,x,\zeta)$ with $\xi_x$ curvilinear,
are easily seen to deform into this open subset.

{\em Claim: $\dim(W\setminus W')<2\ell+m+n+1$. In particular, the complement
of $W'$ in $W$ cannot support any contribution to $w$.}

Indeed, the set $T=\{(\xi,x,x,\zeta)|\xi-\zeta=(n+m)x\}$ has a stratification $T=\coprod_{i\geq0} T_i$, where the stratum $T_i$ is the locally closed
set of all tuples with $\length(\zeta_x)=i$. Let $T'_0\subset T_0$ be the 
closed subset that consists of tuples where $\xi_x$ is not curvilinear. Then
$W\setminus W'\subset T_0'\cup T_1\cup T_2\ldots$.
Now $T_0$ is irreducible of dimension $2\ell+(n+m+1)$, and $T_0'$ is a proper
closed subset and therefore has strictly smaller dimension. The assertion now
follows from Lemma \ref{Z0andZ1}.

{\em Claim: The intersection of $p_{123}^*[Q^{[\ell+m+n]}]$ and $p_{345}^*[Q^{[\ell+m,m]}]$ is transversal at general points of $V'$.}

In fact, the intersection is transversal at all points with $x\neq y$ and $\xi$ curvilinear.

We conclude, that the intersection cycle $v$ equals $[\overline{V'}]+r$, where
$r$ is a cycle supported on $p_{1245}^{-1}(W\setminus W')$ and therefore
irrelevant for our further computations for dimension reasons.
Let us return to the definition of the cycle $w$.

Identifying $V'$ and $W'$ we see that the variety $W'$ parametrises
three families
$$Z\subset \Sigma\subset \Xi\subset W'\times X$$
of subschemes in $X$. In terms of these we can summarise the discussion above
by stating that $\qq_n'\qq_m$ is induced by the cycle
$$c_1(p_*\ki_{\Sigma/\Xi})\cdot[W']\in A_*(W').$$

Having reached this point we pause to reflect what changes in this picture
if we exchange the order of the operators $\qq_n$ and $\qq_m$. Up to the
usual twist $T$ that flips the factors $X$ in $X^{[\ell+m+n],[1],[1],[\ell]}$, 
not a iota is changed in $W'$. Indeed, 
$W'$ parametrises not only three but rather four families of subschemes
$$\begin{array}{ccccc}
&&\Sigma'\\
&\nearrow&&\searrow\\
Z&&&&\Xi\\
&\searrow&&\nearrow\\
&&\Sigma''
\end{array}$$
where $\Sigma'$ and $\Sigma''$ are characterised by the property that at
a point $s=(\Xi_s,x,y,Z_s)\in W'$ the subschemes 
$\Sigma_s', \Sigma_s''\subset\Xi_s$ are the unique ones with 
$$\Sigma_s'-Z_s=mx\,,\quad \Xi_s-\Sigma_s'=ny$$
and 
$$\Sigma_s''-Z_s=ny\,,\quad \Xi_s-\Sigma_s''=mx.$$
This means: the commutator $[\qq_n',\qq_m]$ is induced by the cycle
$$\Big(c_1(p_*\ki_{\Sigma'/\Xi})-c_1(p_*\ki_{Z/\Sigma''})\Big)\cdot[W']
\in A_{2\ell+n+m+1}(X^{[\ell+n+m],[1],[1],[\ell]}).$$

The ideal sheaves corresponding to the various inclusions between the families $Z$,
$\Sigma'$, $\Sigma''$ and $\Xi$ are related by the following commutative
diagram of short exact sequences
$$\begin{array}{ccccccccc}
0&\lra&\ki_{\Sigma'/\Xi}&\lra&\ki_{Z/\Xi}&\lra&\ki_{Z/\Sigma'}&\lra&0\\[1ex]
&&\phi\Big\downarrow\phantom{\phi}&&\Big\|&&\Big\uparrow\\[.5ex]
0&\lla&\ki_{Z/\Sigma''}&\lla&\ki_{Z/\Xi}&\lla&\ki_{\Sigma''/\Xi}&\lla&0\\
\end{array}.$$
The homomorphism
$$p_*\phi:p_*\ki_{\Sigma'/\Xi}\to p_*\ki_{Z/\Sigma''}$$
is an isomorphism off the diagonal $\{x=y\}\subset W'$. On the other hand
the closure of $W'\cap\{x=y\}$ equals the image of the `diagonal' embedding
$Q^{[\ell+m+n,\ell]}\to X^{[\ell+m+n],[1],[1],[\ell]}$.
It follows that
$$\Big(c_1(p_*\ki_{\Sigma'/\Xi})-c_1(p_*\ki_{Z/\Sigma''})\Big)\cdot[W']=-\mu
\cdot[Q^{[\ell+m+n,\ell]}]$$
where $\mu$ is the length of $\mbox{coker}(p_*\phi)$ at the generic point of
the  variety $Q^{[\ell+m+n,\ell]}$.
This proves
$$[\qq_n'(\alpha),\qq_m(\beta)]=-\mu\cdot\qq_{n+m}(\alpha\beta),$$
and it remains to show that
$$\mu=nm.$$

A general point $d=(\xi,x,y,\zeta)$ of $Q^{[\ell+m+n,\ell]}$ is of the form
$(\zeta\cup\eta,x,x,\zeta)$ where $\eta\cap\zeta=\emptyset$ and $\eta$ is
a curvilinear subscheme supported at $x$. As the computation is local in $X$
we may apply the same reduction process as in the proof of Theorem \ref{eqcommutator}: we may assume that $\ell=0$,
that $X=\IA^2=\Spec\IC[z,w]$, $x=(0,0)$ and $I_\zeta=(w,z^n)$. Then there is an 
open neighbourhood of this point $d$ in $W'$ which isomorphic to $\IA^{n+m+2}=\Spec\IC[a_0,\ldots,a_{n+m-1},s,t]$
such that the families $\Xi,\Sigma'$ and $\Sigma''$ are given by the ideals
$$I_\Xi=(w-f(z),(z-t)^n(z-s)^m),\quad I_{\Sigma'}=(w-f(z),(z-s)^m)$$
and
$$I_{\Sigma''}=(w-f(z),(z-t)^n),$$
where $f(z)=a_0+a_1z+\ldots+a_{n+m-1}z^{n+m-1}$. We find
$$p_*\ko_{\Sigma''}=\IC[\underline{a},s,t][z]/(z-t)^n$$
and 
$$p_*\ki_{\Sigma'/\Xi}=(z-s)^m\cdot\IC[\underline{a},s,t][z]/(z-s)^m(z-t)^n.$$
The cokernel of 
$$p_*\phi:(z-s)^m\cdot\IC[\underline{a},s,t][z]/(z-s)^m(z-t)^n
\lra\IC[\underline{a},s,t][z]/(z-t)^n$$
is isomorphic to the $\IC[\underline{a},s,t]$-module
$$\IC[\underline{a},s,t][z]/((z-s)^m,(z-t)^n)\isom
\IC[\underline{a},s+t][z-s,z-t]/((z-s)^m,(z-t)^n).$$ 
This module is supported along the diagonal $\{s=t\}$ (as we expected), and
its stalk at the generic point of the diagonal has length $nm$ (as we had
to prove).

{\em Step 2:} Assume that $m$ is positive and $-m<n<0$. First one shows 
as above that the commutator $[\qq_n',\qq_m]$ is
induced by cycles in $A_{2\ell+n+m+1}(X^{[\ell+m+n],[1],[1],[\ell]})$ for
each $\ell\geq0$, which are supported on the diagonally embedded varieties
$Q^{[\ell+m+n,\ell]}$, so that
$$[\qq_n'(\alpha),\qq_m(\beta)]=-c_{n,m}\cdot \qq_{n+m}(\alpha\beta)$$
for certain constants $c_{n,m}$. In order to determine these constants
we apply the commutator $[\,.\,,\qq_{-n-m}(1)]$.
Then the oscillator relations yield for the right hand side  
$$-c_{n,m}(n+m)\int_X\alpha\beta\cdot\id_\IH.$$
On the other hand
\begin{eqnarray*}[\,[\qq_{n}'(\alpha),\qq_m(\beta)],\qq_{-n-m}(1)]&=&
[\,[\qq_{n}'(\alpha),\qq_{-n-m}(1)],\qq_m(\beta)]
\end{eqnarray*}
Now
\begin{eqnarray*}[\qq_{n}'(\alpha),\qq_{-n-m}(1)]&=&(-1)^{m}[(\qq_{-n}^\dagger)'(\alpha),\qq_{n+m}^\dagger(1)]\\
&=&-(-1)^m[\qq_{n+m}(1),\qq_{-n}'(\alpha)]^\dagger,
\end{eqnarray*}
which by Step 1 equals $(-1)^mn(n+m)\qq_m(\alpha)^\dagger=n(n+m)\qq_{-m}(\alpha)$.
Hence
\begin{eqnarray*}
[\,[\qq_{n}'(\alpha),\qq_m(\beta)],\qq_{-n-m}(1)]&=&n(n+m)[\qq_{-m}(\alpha),
\qq_m(\beta)]\\
&=&n(n+m)(-m)\int_X\alpha\beta\cdot\id_\IH.
\end{eqnarray*}
Choose classes $\alpha,\beta$ with $\int_X\alpha\beta\neq0$. It follows that
$c_{n,m}=nm$.

{\em Step 3:} The general case can now be reduced formally to the cases already treated. The assertion is certainly trivial if either $n=0$ or $m=0$.
If the assertion is known to be true for some pair $(n,m)$, we may apply the
operation $\dagger$ to both sides and find:
\begin{eqnarray*}[\qq_{-n}'(\alpha),\qq_{-m}(\beta)]
&=&(-1)^{n+m}[(\qq_n^\dagger)'(\alpha),\qq^\dagger_m(\beta)]\\
&=&-(-1)^{n+m}[(\qq_n')^\dagger(\alpha),\qq_m^\dagger(\beta)]\\
&=&(-1)^{n+m}[\qq_n'(\alpha),\qq_m(\beta)]^\dagger=-nm\cdot(-1)^{n+m}\qq_{n+m}^\dagger(\alpha\beta)\\
&=&(-n)(-m)\cdot\qq_{-n-m}(\alpha\beta).
\end{eqnarray*}
This and the identity
$$[\qq_n'(\alpha),\qq_m(\beta)]=(-1)^{|\alpha|\cdot|\beta|}[\qq_m'(\beta),\qq_n(\alpha)]$$
allow us to reduce anything to cases checked in Step 1 and Step 2.
\qed

In order to prove part 1 of Theorem \ref{maintheorem}, it remains to treat
the case $n+m=0$. This will be done in two steps. First, we prove a qualitative
statement about the structure of the `correction term', and afterwards we 
determine the precise value of the `coefficient' $K_n$:
 
\begin{proposition}\label{n-plus-m-gleichnull}--- There exist rational divisors $K_n\in Pic(X)\tensor\IQ$, $n\in \IZ$, with $K_0=0$ and $K_{-n}=K_n$ and such that
\beeq{defofKn}[\qq_n'(\alpha),\qq_{-n}(\beta)]=n^2\cdot\int_XK_n\alpha\beta\cdot\id_{\IH}
\eneq
for all $\alpha,\beta\in H^*(X)$.
\end{proposition}

\prf There is nothing to prove for $n=0$. Moreover, $$[\qq_n'(\alpha),\qq_{-n}(\beta)]=(-1)^{|\alpha|\cdot|\beta|}\cdot[\qq_{-n}'
(\beta),\qq_n(\alpha)].$$
It follows that if there is a divisor $K_n$ so that
(\ref{defofKn}) holds for $n$, then (\ref{defofKn}) also holds for $-n$ with
the choice $K_{-n}=K_n$. Hence it suffices to prove the proposition for
positive integers $n$.

Let $\ell$ be a nonnegative integer and consider the diagram
$$\begin{array}{ccccc}
X^{[\ell],[1],[\ell+n]}&\stackrel{p_{123}}{\verylongleftarrow{3em}}&
X^{[\ell],[1],[\ell+n],[1],[\ell]}&
\stackrel{p_{345}}{\verylongrightarrow{3em}} X^{[\ell+n],[1],[\ell]} \\
&&\phantom{\scriptstyle p_{1245}}\Bigg\downarrow
{\scriptstyle p_{1245}}\\[3ex]
&&X^{[\ell],[1],[1],[\ell]}.
\end{array}$$
Let 
$$v:=p_{123}^*[Q^{[\ell,\ell+n]}]\cdot p_{345}^*[Q^{[\ell+n,\ell]}]\in
A_{2\ell+2}(V),$$
{}$$V:=p_{123}^{-1}(Q^{[\ell,\ell+n]})\cap 
p_{345}^{-1}(Q^{[\ell+n,\ell]}).$$
According to Lemma \ref{wasistdieableitung}, the operator $\qq_{-n}'\qq_n$
is induced by the class
$$w=(-1)^np_{1245*}(p_{123}^*c_1(\ki_{\ell,\ell+n})\cdot v)\in A_{2\ell+1}(W), W:=p_{1245}(V).$$
Consider the diagonal part $W\cap\{x=y\}$ first. It is contained in
$\bigcup_{i\geq0}T_i$, where
$T_i=\{(\xi,x,x,\zeta)|\ell(\xi_x)=\ell(\zeta_x)=i\}$.
The closure of $T_0$ is the diagonal $\Delta\isom
X^{[\ell]}\times X\subset X^{[\ell],[1],[1],[\ell]}$ and is therefore 
irreducible of dimension $2\ell+2$. Whereas for $i\geq 1$, the set $T_i$
embeds into the irreducible variety $X^{[\ell-i]}\times (X^{[i]}_0\times_X
X^{[i]}_0)$ of dimension $2(\ell-i)+(i+1)+(i+1)-2=2\ell$.

The off-diagonal part $W\cap\{x\neq y\}$ is empty if $\ell<n$. If $\ell\geq n$
it has precisely one irreducible component $W'$
of maximal dimension $2\ell+2$: it contains as a dense subset the image of
the embedding
$$\{(\eta,\xi',\zeta')\in X^{[\ell-n]}\times X^{[n]}_0\times X^{[n]}_0| \eta,\xi'\mbox{ and }\zeta'\mbox{ are pairwise disjoint}\}\lra W,$$
{}
$$(\sigma,\xi',\zeta')\mapsto(\sigma\cup\xi,\rho(\xi'),\rho(\zeta'),\sigma\cup
\zeta').$$
Since the function $(\xi,x,y,\zeta)\mapsto \ell(\xi_x)$ is semicontinuous and
is at least $n$ on $W'$, it follows that $\overline{W'}\cap \Delta$ is
contained in $\bigcup_{\nu\geq n}T_n$. In particular, this intersection has
dimension $\leq 2\ell$. As we want to compute a cycle of degree $2\ell+1$,
we may restrict our attention to the open part $W'$ and may disregard the complement
of $W'$ in its closure.

$p_{1245}:p_{1245}^{-1}(W')\to W'$ is an isomorphism, which we use to
identify $W'$ and the off-diagonal part of $V$. Now $W'$ parametrises
four flat families of subschemes on $X$: besides the families $\Xi$ and $Z$ 
of fibrewise length $\ell$, these are the families $\Xi\cap Z$ and $\Xi\cup Z$ of
fibrewise length $\ell-n$ and $\ell+n$. The contribution of $W'$ to
$w$ is the class
$$(-1)^nc_1(p_*\ki_{\Xi/\Xi\cup Z})\cdot[W']\in A_{2\ell+1}(W').$$

Reversing the order of the operators $\qq_{-n}'$ and $\qq_n$ shows that
the part of the cycle $u$ inducing the commutator $[\qq_{-n}',\qq_n]$, that is 
supported on $W'$, is the class
$$(-1)^n\Big(c_1(p_*\ki_{\Xi/\Xi\cup Z})-c_1(p_*\ki_{\Xi\cap Z/\Xi})\Big)\cdot[W'].$$
Since the ideal sheaves $\ki_{\Xi/\Xi\cup Z}$ and $\ki_{\Xi\cap Z/\Xi}$ are
isomorphic, this class is zero.

Thus we may fully concentrate on the contribution of the diagonal part $\Delta$.
(Also note that for the reversed order $\qq_n\qq_{-n}'$ any diagonal parts
must be contained in $\bigcup_{\nu\geq n}T_\nu$ and are therefore too small
and irrelevant.)

The complement of the open subset $T_0\isom X^{[\ell]}\times X\setminus\Xi_\ell$
in $\Delta_0$ has codimension $\geq 2$. Locally near $p_{1245}^{-1}(T_0)$
there are isomorphisms between $X^{[\ell+n,\ell]}$ and $X^{[\ell]}\times 
X^{[n]}$, and similarly between $Q^{[\ell+n,\ell]}$ and $X^{[\ell]}\times X^{[n]}_0$. Hence if $\bar w\in A_1(X)$ is the intersection cycle for the 
special case $\ell=0$, then the general cycle is simply given by $w=[X^{[\ell]}]\times \bar w\in A_{2\ell+1}(X^{[\ell]}\times X)$. But that
was all we had to prove: a cycle of this form induces the linear map
$$\alpha\tensor\beta\tensor y\mapsto \int_{\bar w}\alpha\beta \cdot y,\qquad
\alpha,\beta\in H^*(X;\IQ), y\in \IH.$$
\qed

\begin{corollary}--- For all positive integers $n$ one has
$$\qq_n'(\alpha)=n\LL_n(\alpha)+n\qq_n(K_n\alpha).$$
\end{corollary}

\prf Use the same argument as in the first paragraph of the proof of the main
theorem after Corollary \ref{q1machtsallein}.\qed

To finish the proof of Theorem \ref{maintheorem} it remains to show:

\begin{proposition}\label{identifyKn}--- For all positive integers $n$ the rational divisor
defined by Proposition \ref{n-plus-m-gleichnull} is given by
$$K_n=\frac{n-1}{2}K,$$
where $K$ is the canonical class of the surface $X$.
\end{proposition}
 
This will be done in the next section.

\subsection{The vertex operator, completion of the proof}

\begin{definition}--- Let $\gamma\in H^*(X)$ be an element which is of even
degree though not necessarily homogeneous, and let $t$ be a formal parameter.
Define operators $S_m(\gamma)$, $m\geq0$, by
$$S(\gamma,t):=\sum_{m\geq0}S_m(\gamma)t^m:=
\exp\left(\sum_{n>0}\frac{(-1)^{n-1}}{n}\qq_n(\gamma)\cdot t^n\right).$$
\end{definition}

Since $\gamma$ is of even degree by assumption, any two operators 
$\qq_n(\gamma)$ and $\qq_{n'}(\gamma)$ commute in the ordinary, i.e.\
`ungraded' sense. In particular, there is no ambiguity in the meaning of the
expression on the right hand side in the definition.

The geometric meaning of the operators $S_m$ is explained by the following theorem: let $C$ be a smooth curve in $X$. There is an
induced closed embedding $S^nC=C^{[n]}\to X^{[n]}$. Let $[C]\in H^*(X)$ and
$[C^{[n]}]\in H^*(X^{[n]})$ be the corresponding cohomology classes, i.e., the
Poincar\'e dual classes of the fundamental classes of these varieties.

\begin{theorem}[Nakajima, Grojnowski]\label{Nakonsymmprod}--- The following relation holds for all nonnegative integers $n$:
$$[C^{[n]}]=S_n([C])\cdot \vacuum.$$
\end{theorem}

For proofs see \cite{NakajimaLectures} and \cite{Grojnowski}.\qed

\begin{lemma}--- Let $\gamma\in H^*(X)$ be an element of even degree.
Then
$$S'(\gamma,t)=S(\gamma,t)\cdot\sum_{n>0}(-1)^{n-1}
t^n\left\{\LL_n(\gamma)+\qq_n\Bigl(\gamma K_n+\gamma^2\frac{n-1}{2}\Bigr)
\right\}.$$
\end{lemma}

\prf Assume first that $\gotha$ is an operator of even degree, and that
$[\gotha',\gotha]$ commutes with $\gotha$. Then
\begin{eqnarray*}
\left(\sum_{n=0}^\infty\frac{\gotha^n}{n!}\right)'&=&\sum_{n=1}^\infty\frac{1}{n!}
\sum_{i=1}^n \gotha^{i-1}\cdot\gotha'\cdot\gotha^{n-i}\\
&=&\sum_{n=1}^\infty\frac{1}{n!}\cdot\left\{n\gotha^{n-1}\gotha'+
\sum_{i=1}^n\gotha^{n-2}\cdot(n-i)\cdot[\gotha',\gotha]\right\}\\
&=&\sum_{n=0}^\infty\frac{\gotha^n}{n!}\cdot\gotha'+\sum_{n=1}^\infty\frac{
\gotha^{n-2}}{n!}\binom{n}{2}[\gotha',\gotha]\\
&=&\exp(\gotha)\cdot\left\{\gotha'+\frac{1}{2}[\gotha',\gotha]\right\}.
\end{eqnarray*}

Next, let $\gotha_\nu$ be a family of commuting operators of even degree such
that any $[\gotha'_\nu,\gotha_\mu]$ commutes with every $\gotha_\xi$. Then it 
follows from Step 1 and $$[\gotha'_\mu,\exp(\gotha_\nu)]=\exp(\gotha_\nu)\cdot[\gotha'_\mu,\gotha_\nu]$$
that
$$\left(\exp\Big(\sum_\nu\gotha_\nu\Big)\right)'=
\exp\Big(\sum_\nu\gotha_\nu\Big)
\cdot\left\{\sum_{\nu}\gotha_\nu'+
\frac{1}{2}\sum_{\nu,\mu}[\gotha_\nu',\gotha_\mu]\right\}.$$
Now apply this formula to the family $\gotha_\nu=\frac{(-1)^{\nu-1}}{\nu}\qq_{\nu}(\gamma)t^\nu$
and use our previous results
$\gotha'_\nu=(-1)^{\nu-1}t^\nu(\LL_n(\gamma)+\qq_\nu(K_\nu\gamma))$
and $[\gotha'_\nu,\gotha_\mu]=-(-t)^{\nu+\mu}\qq_{\nu+\mu}(\gamma^2)$. 
One gets $S'(\gamma,t)=S(\gamma,t)\cdot(*)$ with
\begin{eqnarray*}
(*)&=&\sum_{n>0}(-1)^{n-1}t^n\big(\LL_n(\gamma)+\qq_n(K_n\gamma)\big)
-\frac{1}{2}\sum_{\nu,\mu>0}(-t)^{\nu+\mu}\qq_{\nu+\mu}(\gamma^2) \\
&=&\sum_{n>0}(-1)^{n-1}t^{n}\cdot
\left\{\LL_n(\gamma)+\qq_n(K_n\gamma+\frac{1}{2}N_n\gamma^2)\right\} \end{eqnarray*}
where $N_n$ is the number of pairs of positive integers $\nu$ and $\mu$ that
add up to $n$, i.e., $N_n=n-1$. 
\qed

Let $C\subset X$ be a smooth projective curve. The boundary $\partial X^{[n]}$
intersects $C^{[n]}$ generically transversely in the boundary 
$\partial C^{[n]}$ of $C^{[n]}$, i.e.\ in the set of all tuples with multiple
points. The subvarieties $X^{[n]}_0$ and $\partial C^{[n]}$ have complementary
dimensions $n+1$ and $n-1$ in $X^{[n]}$ and we may compute the intersection number
$$I:=\int_{X^{[n]}}[X^{[n]}_0]\cup [\partial C^{[n]}].$$
We will do this first using our algorithmic language, and afterwards using a 
geometric argument. The comparison of the two results will lead to the 
identification of the divisors $K_n$.

\begin{lemma}\label{ersterWega}--- $[X^{[n]}_0]=\qq_n(1_X)\cdot\vacuum$ ~~and~~
$[\partial C^{[n]}]=-2\cdot S'_n([C])\cdot\vacuum$.
\end{lemma}

\prf The first assertion follows from the definition of the operators
$\qq_n$.
By Nakajima's Theorem, $S_n([C])\cdot\vacuum$ is the
class of the submanifold $C^{[n]}\subset X^{[n]}$, and hence according to Lemma
\ref{DnX}:
$$S_n'([C])\cdot\vacuum=\dd\cdot S_n([C])\cdot\vacuum=
-\frac{1}{2}[\partial X^{[n]}]\cdot[C^{[n]}]=-\frac{1}{2}[\partial C^{[n]}].$$
\qed

\begin{lemma}\label{ersterWegb}---
$$\int_{X^{[n]}}(\qq_n(1_X)\cdot\vacuum)\cdot(S_n'([C])\cdot\vacuum)=\int_X \left\{nK_nC+\binom{n}{2}C^2\right\}.$$
\end{lemma}

\prf Indeed,
\begin{eqnarray*}
\int_{X^{[n]}}(\qq_n(1_X)\cdot\vacuum)\cdot(S_n'([C])\cdot\vacuum)
&=&(-1)^n\int_{X^{[0]}}\qq_{-n}(1_X)S_n'([C])\cdot\vacuum\\
&=&(-1)^n\int_{X^{[0]}}[\,\qq_{-n}(1_X),S'_n([C])\,]\cdot\vacuum,
\end{eqnarray*}
since $\qq_{-n}(1_X)\cdot\vacuum=0$.
Now $\qq_{-n}$ commutes with any product $\qq_{i_1}\cdot\ldots
\cdot\qq_{i_s}$ if $s\geq2$, $i_j>0$ and $\sum_ji_j=n$. Thus the only summand
in $S'_n$ that contributes to the commutator with $\qq_{-n}$ is $(-1)^{n-1}\qq_n
(C(K_n+C(n-1)/2))$. Hence
$$[\qq_{-n}(1_X), S'_n([C])]
=(-1)^nn\int_XC\left(K_n+\frac{n-1}{2}C\right)\cdot\id_{\IH}$$
This proves the lemma.\qed

Next, we give the geometric computation of $I$:
 
\begin{lemma}\label{zweiterWeg}--- 
$$\int_{X^{[n]}}[X_0^{[n]}]\cdot [\partial C^{[n]}]=-n(n-1)\cdot C(C+K).$$
\end{lemma}

\prf We have $[X^{[n]}_0]\cdot[\partial C^{[n]}]=[\partial X^{[n]}]\cdot(
[X^{[n]}_0]\cdot[C^{[n]}])$.
The intersection of $X^{[n]}$ and $C^{[n]}$ is transversal and is equal to 
the image of the closed immersion $\Delta:C\to C^{[n]}$ sending a point $c$ to 
the unique subscheme of $C$ of length $n$ that is supported in $c$.
Thus
$$I=\deg(\ko_{X^{[n]}}(\partial X^{[n]})|_{\Delta(C)}=\deg(\ko_{C^{[n]}}
(\partial C^{[n]})|_{\Delta(C)}.$$
The embedding $\Delta$ factors through the diagonal embedding
$C\to C^n$ and the quotient map $\pi:C^n\to C^{[n]}$. Moreover, if
$\pr_{ij}:C^n\to C^2$ denotes the projection to the product of the $i$-th and 
$j$-th factor,
$$\pi^*(\ko_{C^{[n]}}(\partial C^{[n]}))\isom
\left(\bigotimes_{i<j}^n pr_{ij}^*\ko_{C\times C}(\Delta_C)\right)^{\otimes2}.$$
From this we conclude:
\begin{eqnarray*}I=\deg(\Delta^*\ko_{C^{[n]}}(\partial C^{[n]}))&
=&2\cdot\binom{n}{2}\deg(\ko_{C\times C}(\Delta_C)|_{\Delta_C}))\\
&=&-n(n-1)\cdot C(C+K).
\end{eqnarray*}
\qed

{\kursiv Proof of Proposition \ref{identifyKn}}. From Lemma \ref{ersterWega}
and Lemma \ref{ersterWegb} we conclude
$$I=(-2)\cdot C(nK_n+\binom{n}{2}C ).$$
Comparison with Lemma \ref{zweiterWeg} shows that $K_n=\frac{n-1}{2}K$.\qed

This finishes the proof of Theorem \ref{maintheorem}.

\section{Towards the ring structure of $\IH$}\label{RingStructure}


\subsection{Tautological sheaves}
 
  
There is a natural way to associate to a given vector bundle on $X$ a
series of tautological' vector bundles on the Hilbert schemes $X^{[n]}$, 
$n\geq0$. The Chern classes of the tautological bundles may 
be grouped together to form operators on $\IH$.

Consider the standard diagram
$$\begin{array}{ccccc}
\Xi_n&\subset&X^{[n]}\times X&\stackrel{q}{\lra}&X\\
&&{\scriptstyle p}\Big\downarrow\phantom{\scriptstyle p}\\[1ex]
&&X^{[n]}
\end{array}
$$
Let $F$ be a locally free sheaf on $X$. For each $n\geq 0$ the associated
{\em tautological bundle} on $X^{[n]}$ is defined as
$$F^{[n]}:=p_*(\ko_{\Xi_n}\tensor q^*F).$$
Since $p$ is a flat finite morphism of degree $n$,
$F^{[n]}$ is locally free with
$$\rk(F^{[n]})=n\cdot \rk(F).$$
Note that $F^{[0]}=0$ and $F^{[1]}=F$.

Furthermore, if $\ses{F_1}{F}{F_2}$ is a short exact sequence of locally free
sheaves on $X$, the corresponding sequence
$\ses{F_1^{[n]}}{F^{[n]}}{F_2^{[n]}}$ is again exact. 
Hence sending the class $[F]$ of a locally free sheaf $F$ to $[F^{[n]}]$
gives a group homomorphism
$$-^{[n]}:K(X)\lra K(X^{[n]}).$$

\begin{definition}--- Let $u$ be a class in $K(X)$. Define operators
$$\cc(u)\in \End(\IH)\quad\mbox{ and }\quad\ch(u)\in \End(\IH)$$
as follows: For each $n\geq0$, the action on $H^*(X^{[n]};\IQ)$ is given by 
multiplication with the total Chern class $c(u^{[n]})$ and the 
Chern character $ch(u^{[n]})$, respectively. 
\end{definition}

Let
$$\cc(u)=\sum_{k\geq0}\cc_k(u)\quad\mbox{ and }\quad\ch(u)=
\sum_{k\geq 0}\ch_k(u)$$
be the decompositions into homogeneous components of bidegree $(0,2k)$.
Since all of these operators are of even degree and only act `vertically'
on $\IH$ by multiplication, they commute with each other and in particular 
with the previously defined boundary operator $\dd=\cc_1(\ko_X)$.

Moreover, we have
$$\cc(u+v)=\cc(u)\cdot\cc(v)\quad\mbox{ and }\quad\ch(u+v)=\ch(u)+\ch(v)$$
for all $u,v\in K(X)$.

\begin{theorem}\label{cqcommutator}--- Let $u$ be a class in $K(X)$
of rank $r$ and let $\alpha\in H^*(X)$. Then
$$[\ch(u),\qq_1(\alpha)]=\exp(\ad\dd)(\qq_1(ch(u)\alpha)),$$
or, more explicitly,
$$[\ch_n(u),\qq_1(\alpha)]=\sum_{\nu=0}^n\frac{1}{\nu!}\qq_1^{(\nu)}
(ch_{n-\nu}(u)\alpha).$$
Similarly,
$$\cc(u)\cdot\qq_1(\alpha)\cdot\cc(u)^{-1}=\sum_{\nu,k\geq0}\binom{r-k}{\nu}\qq_1^{(\nu)}
(c_k(u)\alpha).$$
\end{theorem}

\prf We may assume that $u$ is the class of a locally free sheaf $F$.
Recall the standard diagram for the incidence variety $X^{[\ell,\ell+1]}$:
$$
\begin{array}{ccccc}
X&\stackrel{\rho}\longleftarrow&X^{[\ell,\ell+1]}&\stackrel{\psi}{\longrightarrow}&X^{[\ell+1]}\\
&&{\scriptstyle \phi}\Big\downarrow\phantom{\scriptstyle \phi}\\
&&X^{[\ell]}
\end{array}
$$
The variety $X^{[\ell,\ell+1]}$ parametrises two families of subschemes of $X$. Their structure sheaves fit into an exact sequence
$$\ses{\rho_X^*\ko_{\Delta_X}\tensor p^*\ko_{X^{[\ell,\ell+1]}}(-E)}
{\psi_X^*(\ko_{\Xi_{\ell+1}})}{\phi_X^*(\ko_{\Xi_{\ell}})},$$
where $p:X^{[\ell,\ell+1]}\times X\to X^{[\ell,\ell+1]}$ is the projection
and $E$ is the exceptional divisor. 
Applying the functor $p_*(\,\cdot\,\tensor q^*F)$ to this exact sequence yields
\beeq{Geirssequence}
\ses{\rho^*F\tensor\ko_{X^{[\ell,\ell+1]}}(-E)}{\psi^*F^{[\ell+1]}}
{\phi^*F^{[\ell]}}.
\eneq
Let $\lambda=c_1(\ko_{X^{[\ell,\ell+1]}}(-E))$. Then
$$\psi^*ch(F^{[\ell+1]})=\phi^*ch(F^{[\ell]})+\rho^*ch(F)\cdot \exp(\lambda)$$
and
$$\psi^*c(F^{[\ell+1]})=\phi^*c(F^{[\ell]})\cdot \sum_{\nu,k\geq0}\binom{r-k}{\nu}\lambda^\nu \rho^*c_k(F).$$
It follows for any $x\in H^*(X^{[\ell]};\IQ)$:
\begin{eqnarray*}
\ch(F)\qq_1(\alpha)(x)&=&ch(F^{[\ell+1]})\cdot 
PD^{-1}\psi_*([X^{[\ell,\ell+1]}]\cap \rho^*(\alpha)\phi^*(x))\\
&=&PD^{-1}\psi_*([X^{[\ell,\ell+1]}]\cap\psi^*(ch(F^{[\ell+1]}))\rho^*(\alpha)
\phi^*(x))\\
&=&PD^{-1}\psi_*([X^{[\ell,\ell+1]}]\cap\rho^*(\alpha)\phi^*(ch(F^{[\ell]})x))\\
&&+\sum_{\nu\geq0}\frac{1}{\nu!}PD^{-1}\psi_*(\lambda^\nu
\cdot[X^{[\ell,\ell+1]}]\cap \rho^*(ch(F)\alpha)\phi^*(x))\\
&=&\qq_1(\alpha)(\ch(F)x)+\sum_{\nu\geq0}\frac{1}{\nu!}\qq^{(\nu)}
(ch(F)\alpha)(x).
\end{eqnarray*}
Here we used Lemma \ref{wasistdieableitung} which says that the cycle 
$\lambda^\nu\cdot[X^{[\ell,\ell+1]}]$ 
induces the operator $\qq_1^{(\nu)}$. This is the equation for the Chern 
character. The equation for the total Chern class is proved analogously.
\qed

\begin{corollary}\label{Chernclasses}--- For any $u\in K(X)$ let $\gothC(u)$ be the operator
$$\gothC(u)=\gothc(u)\cdot\qq_1(1_X)\cdot\gothc(u)^{-1}=
\sum_{\nu,k\geq0}\binom{\rk(u)-k}{\nu}\qq_1^{(\nu)}(c_k(u)\alpha).$$
Then
$$\sum_{n\geq0}c(u^{[n]})=\exp(\gothC(u))\cdot\vacuum.$$
\end{corollary}

Note that the right hand side can be explicitly expressed in terms of the basic
operators $\qq_n$ by applying Theorem \ref{maintheorem}.

\prf We have
\begin{eqnarray*}
\sum_{n\geq0}c(u^{[n]})&=&\gothc(u)\sum_{n\geq0}1_{X^{[n]}}\\
&=&\gothc(u)\exp(\qq_1(1_X))\cdot\vacuum\\
&=&\gothc(u)\exp(\qq_1(1_X))\gothc(u)^{-1}\cdot\vacuum\\
&=&\exp(\gothc(u)\qq_1(1_X))\gothc(u)^{-1})\cdot\vacuum\\
&=&\exp(\gothC(u))\cdot\vacuum.
\end{eqnarray*}
\qed

\begin{remark}--- The sequence (\ref{Geirssequence}) was used by Ellingsrud 
in a recursive method to determine Chern classes and Segre classes
of tautological bundles (unpublished, but see \cite{TikhomirovStandardBundles},\cite{EllingsrudGoettscheLehn}). He expresses the classes $(\phi,\rho)_*c(E)$ in terms of the 
Segre classes of the universal family $\Xi_{[n]}\subset X\times X^{[n]}$.
Thus one needs to control the behaviour of these Segre classes under the 
induction procedure. This method yields qualitative results on the {\em
structure} of certain classes and integrals, but all attempts to get numbers
have ended so far in unsurmountable combinatorial difficulties.\qed
\end{remark}
 
\begin{remark}--- The results of the present and the previous section provide
an algorithmic description of the multiplicative action of the subalgebra $\ka\subset\IH$ which is generated by the Chern classes of all tautological bundles:
The elements $\qq_{i_1}(\alpha_1)\cdot\ldots\qq_{i_s}(\alpha_s)\cdot\vacuum$
generate $\IH$ as a $\IQ$-vector space. By Corollary \ref{q1machtsallein}, 
each such
element can be written as a linear combination of expression $w\cdot\vacuum$,
where $w$ is a word in an alphabet consisting of $\dd$ and operators $\qq_1(\alpha)$, $\alpha\in H^*(X;\IQ)$. By Theorem \ref{cqcommutator} 
the commutator of $\ch(F)$ with any of these is again a word in this alphabet.
And finally, Theorem \ref{maintheorem} shows how such a word can be expressed
in terms of the basic operators $\qq_n$. Admittedly, without a further 
understanding of the algebraic structure this description is useful for
computations in $H^*(X^{[\ell]};\IQ)$ only for small values of $\ell$ or if one
implements it in some computer algebra system. The following sections 
deal with special situations where one can say more.
\end{remark}


\subsection{The line bundle case}


The results of the previous section suffice to compute the Chern classes
of the tautological bundles $L^{[n]}$ associated to a line bundle $L$ in
terms of the basic operators.

\begin{theorem}\label{Lotharconjecture}--- Let $L$ be a line bundle on $X$. Then
$$\sum_{n\geq0}c(L^{[n]})=\exp\left(\sum_{m\geq1}
\frac{(-1)^{m-1}}{m}\qq_m(c(L))\right)\cdot\vacuum.$$
\end{theorem}

\begin{remark}--- Expanding the term on the right hand side, one realises that
the cohomological degree of any summand contained in $H^*(X^{[n]};\IQ)$ is
$\leq 2n$, and, moreover, the maximal degree
$2n$ can only be attained if the arguments of all operators $\qq_\nu$ involved
have degree 2. In other words, considering elements of
top degree only, the equation of the theorem specialises to
\beeq{Lotharconjecturespecialcase}
\sum_{n\geq0}c_n(L^{[n]})=\exp\left(\sum_{m\geq1}\frac{(-1)^{m-1}}{m}\qq_m(c_1(L))\right)\cdot\vacuum.
\eneq
This is Nakajima's result \ref{Nakonsymmprod}: for suppose $C\subset X$ is a
smooth curve and $L=\ko_X(C)$. If $\xi\in X^{[n]}$, the natural homomorphism
$\ko_X\to\ko_{\xi}(C)$ vanishes if and only if $\xi\subset C$. Hence the 
vanishing locus of the global vector bundle homomorphism
$$\ko_{X^{[n]}}\lra (\ko_X(C))^{[n]}=L^{[n]}$$
is the subvariety $C^{[n]}$. Therefore $[C^{[n]}]=c_n(L^{[n]})$. Inserting this
into (\ref{Lotharconjecturespecialcase}), we recover Nakajima's
formula \ref{Nakonsymmprod}
$$\sum_{n\geq0}[C^{[n]}]=\exp\left(\sum_{m\geq1}\frac{(-1)^{m-1}}{m}\qq_m([C])
\right)\cdot\vacuum$$
Based on this observation, the theorem was conjectured by L.\ G\"ottsche in
a letter to G.\ Ellingsrud and the author.
\end{remark}

\prfofthetheorem We shall give two variants of the proof which differ slightly 
in flavour. We have seen that the left hand side in the theorem equals
$\exp(\gothC(L))\cdot\vacuum$, where in this case because of $r=1$ we have $$\gothC(L)=\qq_1(c(L))+\qq_1'(1_X).$$

{\em Variant 1}. Expanding the right hand side of
$$\exp(\gothC(L))\cdot\vacuum=\sum_{n\geq0}\frac{1}{n!}(\qq_1(c(L))+\qq_1'(1_X))^n\cdot\vacuum
$$
yields summands which are words in the two symbols
$\qq_1(c(L))$ and $\qq_1'(1_X)$. Moving all factors $\qq_1'(1_X)$ within a
given word as far to the right as possible using the commutation relations of 
the main theorem we can write 
$$\sum_{n\geq0}\frac{1}{n!}(\qq_1(c(L))+\qq_1'(1_X))^n\cdot\vacuum
=\gothA\cdot\vacuum+\gothB\cdot\qq_1'(1_X)\cdot\vacuum=\gothA\cdot\vacuum,$$
where $\gothA$ is a sum of expressions of the form $$
\nu_1!\cdot\ldots\cdot\nu_s!\cdot
\frac{(-1)^{\nu_1-1}\qq_{\nu_1}(c(L))}{\nu_1}\cdots
\frac{(-1)^{\nu_s-1}\qq_{\nu_s}(c(L))}{\nu_s}.$$
Let $\alpha=(1^{\alpha_1}2^{\alpha_2}3^{\alpha_3}\ldots)$ denote a partition and let $|\alpha|:=\sum_{i\geq1}i\alpha_i$, and $\alpha!:=\prod_{i}(i!)^{\alpha_i}$. We get
\beeq{Nalpha}
\sum_{n\geq0}\frac{1}{n!}(\qq_1(c(L))+\qq_1'(1_X))^n\cdot\vacuum
=\sum_{\alpha} N_\alpha\frac{\alpha!}{|\alpha|!}
\prod_{i\geq1}\left(\frac{(-1)^{i-1}\qq_i(c(L))}{i}\right)^{\alpha_i}\cdot\vacuum,
\eneq
where the natural number $N_\alpha$ counts how often the operator
$$
{\alpha!}\prod_{i\geq1}
\left(\frac{(-1)^{i-1}\qq_i(c(L))}{i}\right)^{\alpha_i}
$$
arises from a word in $\qq_1'(1_X)$ and $\qq_1(c(L))$ of length $|\alpha|$.
It is not difficult to see that $N_\alpha$ equals the number of possibilities
to partition a set of $|\alpha|$ elements into subsets in such a way
that there are
$\alpha_i$ subsets of cardinality $i$. Hence
$$N_\alpha:=\frac{1}{\alpha_1!\alpha_2!\cdots}\cdot\frac{|\alpha|!}{\alpha!}.$$
Inserting this into equation (\ref{Nalpha}) above one gets
\begin{eqnarray*}
\sum_{n\geq0}\frac{1}{n!}(\qq_1(c(L))+\qq_1'(1_X))^n\cdot\vacuum
&=&\sum_{\alpha}\prod_{i\geq1}\frac{1}{\alpha_i!}
\left(\frac{(-1)^{i-1}\qq_i(c(L))}{i}\right)^{\alpha_i}\cdot\vacuum\\
&=&\prod_{i\geq1}\sum_{\alpha_i\geq0}\frac{1}{\alpha_i!}
\left(\frac{(-1)^{i-1}\qq_i(c(L))}{i}\right)^{\alpha_i}\cdot\vacuum\\
&=&\prod_{i\geq1}\exp\left(\frac{(-1)^{i-1}\qq_i(c(L))}{i}\right)\cdot\vacuum\\
&=&\exp\left(\sum_{i\geq1}\frac{(-1)^{i-1}}{i}\qq_i(c(L))\right)\cdot\vacuum
.
\end{eqnarray*}
In fact, being a little more careful, one gets
$$\exp(\gothC(L))=\exp\left(\sum_{i\geq1}\frac{(-1)^{i-1}}{i}\qq_i(c(L))\right)
\cdot\exp(\qq_1'(1_X)).$$

{\em Variant 2}. Starting again from the sequence
$$\cc(L)\cdot\qq_1(1_X)=\gothC(L)\cdot\cc(L),$$
we multiply by $\frac{1}{n!}\qq_1(1_X)^nt^n$ from the right and sum up over all $n\geq0$:
\begin{eqnarray*}
\frac{d}{dt}\left(\cc(L)\cdot\sum_{n\geq0}\frac{1}{n!}\qq_1(1_X)^n t^n
\right)\cdot\vacuum=
\cc(L)\cdot\sum_{n\geq0}\frac{1}{n!}\qq_1(1_X)^{n+1}t^n\cdot\vacuum\qquad\\
=\gothC(L)\cdot\left(\cc(L)\cdot
\sum_{n\geq0}\frac{1}{n!}\qq_1(1_X)^n t^n\right)\cdot\vacuum.
\end{eqnarray*}
This means that the series
$$\sum_{n\geq0}c(L^{[n]})t^n=\cc(L)\cdot\exp(\qq_1(1_X)t)\cdot\vacuum$$
satisfies the linear differential equation
\beeq{Lotharconjecturediffequation}
\frac{d}{dt}\gothX=\gothC(L)\cdot\gothX
\eneq
with initial condition
\beeq{Lotharconjectureinitial}\quad\gothX(0)=\vacuum.\eneq
On the other hand, consider the operator
$$S(c(L),t)=\exp\left(\sum_{m\geq1}\frac{(-1)^{m-1}}{m}\qq_m(c(L))t^m\right).$$
We find
$$\frac{d}{dt}S(c(L),t)=S(c(L),t)\cdot \left(\sum_{m\geq0}(-1)^{m}\qq_{m+1}(c(L))t^m\right),$$
and
\begin{eqnarray*}
\lefteqn{\Big[\{\qq_1(1_X+c_1(L))+\qq'_1(1_X)\},S(c(L),t)\Big]}\hspace{4em}\\
&=&S(c(L),t)\cdot\left(\sum_{m\geq1}\frac{(-1)^{m-1}}{m}\Big[
\qq_1'(1_X),\qq_m(c(L))\Big]t^m\right)\\
&=&S(c(L),t)\cdot\left(\sum_{m\geq1}(-1)^m\qq_{m+1}(c(L))t^m\right).
\end{eqnarray*}
This shows
\begin{eqnarray*}
\lefteqn{\{\qq_1(1_X+c_1(L))+\qq'_1(1_X)\}\cdot S(c(L),t)\cdot\vacuum}\hspace{1em}\\
&=&S(c(L),t)\cdot\left(\sum_{m\geq1}(-1)^m\qq_{m+1}(c(L))t^m\right)\cdot\vacuum
\\
&&+S(c(L),t)\cdot\qq_1(c(L))\cdot\vacuum\\
&=&S(c(L),t)\cdot\left(\sum_{m\geq0}(-1)^m\qq_{m+1}(c(L))t^m\right)\cdot\vacuum
\end{eqnarray*}
Hence $S(c(L),t)\cdot\vacuum$ satisfies the system
(\ref{Lotharconjecturediffequation}) and
(\ref{Lotharconjectureinitial}) as well and therefore equals $\cc(L)\cdot\exp(\qq_1(1_X)t)\cdot\vacuum$.
This proves the theorem.\qed


\subsection{Top Segre classes}

The following problem was posed by Donaldson in connection with the computation
of instanton invariants: let $n$ be an integer $\geq 1$, and consider a 
linear system $|H|$ of dimension $3n-2$ inducing a map $X\rationalmap
\IP^{3n-2}$. A zero-dimensional subscheme
$\xi\in X^{[n]}$ does not impose  independent conditions on the linear system
$|H|$ if the natural homomorphism
$$H^0(\IP^{3n-2},\ko_\IP(1))\lra H^0(\xi,\ko_\xi(H))$$
fails to be surjective. The subscheme of all such $\xi\in X^{[n]}$ has virtual
dimension zero, and its class is given by $c_{2n}(W\dual)$,
where $W$ is the virtual vector bundle
$$H^0(\IP^{3n-2},\ko_\IP(H))\tensor\ko_{X^{[n]}}-\ko(H)^{[n]}.$$
Thus the number of those $\xi$ that impose dependent conditions is given by
$$N_n:=\int_{X^{[n]}}c_{2n}(-\ko(H)^{[n]})=\int_{X^{[n]}}\cc(-\ko(H))\cdot
\frac{\qq_1(1_X)^n}{n!}\cdot\vacuum.$$
More explicitly, $N_1$ is the degree of the linear system, $N_2$ is the number
of double points, $N_3$ is the number of trisecants to a surface in $\IP^7$
and $N_4$ is the number of quadruples of points on a surface in $\IP^{10}$
that span a plane.

Problem:
{\em Express $N_n$ in terms of intrinsic invariants of $X$ such as the
degree $d:=H.H$, the intersection $\pi:=H.K$ and $\kappa:=K.K$ and
the topological Euler characteristic $e=c_2(X)$.}

Note that even the fact that such an expression in terms of the given invariants
exists is not evident {\em a priori}. This has been proved by Tikhomirov \cite{TikhomirovStandardBundles}.
It also follows immediately from our approach.

Using our algorithm, we can attack this problem as follows.
Theorem \ref{cqcommutator} yields for $F=-\ko(H)$ and $r=-1$ the formula:
\begin{eqnarray*}
\gothC(-\ko(H))
&=&\sum_{\nu,k\geq0}\binom{-1-k}{\nu}\qq_1^{(\nu)}(c_k(-H))\\
&=&\sum_{\nu\geq0}(-1)^\nu\qq_1^{(\nu)}\left(\sum_{k=0}^{2}
\binom{\nu+k}{k}(-H)^k\right)\\
&=&\sum_{\nu\geq0}(-1)^\nu\qq_1^{(\nu)}((1-H+H^2)^{\nu+1}).
\end{eqnarray*}

It follows as in the proof of Theorem \ref{Lotharconjecture}
that $\cc(-\ko(H))\cdot\exp(\qq_1(1_X)t)\cdot\vacuum$ satisfies the
following differential equation and initial value condition:
$$\frac{d}{dt}\gothX=\gothC(-\ko(H))\gothX\quad\mbox{ and }\quad\gothX(0)=\vacuum.$$
As long as no explicit generating function is available we must be content with
the following semi-explicit solution to the problem:
$$N_n=\frac{1}{n!}\int_{X^{[n]}}\gothC(-\ko(H))^n\cdot\vacuum.$$

\begin{example}---
As a special case, let us compute $N_2$. This is the number of secant lines to 
an embedded surface in $\IP^5$ that pass through a fixed but general point 
$x\in\IP^5$. Hence we should find Severi's double point formula \cite{Severi} (see also \cite{Catanese}). Let $\alpha=1-H+H^2$. Then 
$$2\cdot N_2=\int_{X^{[2]}}\gothC(-\ko(H))^2\cdot\vacuum\quad\mbox{with}\quad
\gothC(-\ko(H))=\sum_{n\geq0}(-1)^n\qq^{(n)}_1(\alpha^{\nu+1}).$$
Since $\qq_1^{(n)}\cdot\vacuum=0$ for all $n>0$ and for all parameters, 
we have $\gothC(-\ko(H))\cdot\vacuum=\qq_1(\alpha)\cdot\vacuum$. Moreover, for degree reasons the infinite sum reduces to
$$\gothC(-\ko(H))^2\cdot\vacuum=(\qq_1(\alpha)-\qq_1'(\alpha^2)+\qq_1''(\alpha^3)-\qq_1'''(\alpha^4)+\qq_1''''(\alpha^5))\qq_1(\alpha)\cdot\vacuum.$$
Using $\qq_1^{(\nu)}(x)\qq_1(y)\cdot\vacuum=-\qq_2^{(\nu-1)}(xy)\cdot\vacuum$ 
this becomes
$$\gothC(-\ko(H))^2\cdot\vacuum=(\qq_1(\alpha)\qq_1(\alpha)+\qq_2(\alpha^3)-\qq_2'(\alpha^4)+\qq_2''(\alpha^5)+\qq_2'''(\alpha^6))\cdot\vacuum.$$
For the higher derivatives $\qq_2^{(n)}$, $n\geq2$, there is the 
following recursion formula: 
\begin{eqnarray*}
\qq_2^{(n)}(x)\cdot\vacuum&=&\big(\qq_1^2(\delta(x))+\qq_2(Kx)\big)^{(n-1)}\cdot\vacuum\\
&=&\big(-\qq_2^{(n-2)}(c_2(X)x)+\qq_2^{(n-1)}(Kx)
\big)\cdot\vacuum.
\end{eqnarray*}
(Recall that the composite map 
$H^*(X)\stackrel{\delta}{\lra}H^*(X)\tensor H^*(X)\stackrel{\cup}{\lra} H^*(X)$ 
is the multiplication with the self intersection of the diagonal, i.e.\ the 
second Chern class $c_2(X)$ of $X$.) Using this formula repeatedly and
keeping in mind that $K.e=K^3=e^2=0$ and $K^2.\alpha^\nu=K^2$, $e.\alpha^\nu=e$,
we finally arrive at
$$\gothC(-\ko(H))^2\cdot\vacuum=(\qq_1(\alpha)^2+\qq_1^2\delta(-\alpha^4+K\alpha^5-K^2+e)
+\qq_2(\alpha^3-K\alpha^4+K^2-e))\cdot\vacuum.$$
Only the first two summands contribute to the integral. Hence
\begin{eqnarray*}
2\cdot\int_{X^{[2]}}\gothC(-\ko(H))^2\cdot\vacuum&=&\left(\int_X\alpha\right)^2-\int_X(\alpha^4-K\alpha^4+K^2-e)\\
&=&d^2-10d-5\pi-\kappa+e.
\end{eqnarray*}
\qed
\end{example}

For higher $n$, the practical calculation of $N_n$
quickly becomes rather difficult. Already the case of $N_3$
surpassed my personal calculation skills. Using MAPLE, I computed $N_n$ for $n\leq 7$. One obtains for example:
\begin{eqnarray*}3!\cdot N_3&=&d^3-30d^2+224d-3d(5\pi+\kappa-e)\\
&&+192\pi+56\kappa-40e,
\end{eqnarray*}
\begin{eqnarray*}
4!\cdot N_4&=&d^4-60d^3+d^2(1196-30\pi+6e-6\kappa)\\
&&-d(7920-1068\pi+220e-284\kappa)+3e^2+1944e-6e\kappa\\
&&-30e\pi+75\pi^2+3\kappa^2+30\kappa\pi-9042\pi
-3300\kappa,
\end{eqnarray*}
\begin{eqnarray*}5!\cdot N_5&=&
d^5-100d^4+d^3(3740+10e-50\pi-10\kappa)\\
&&-d^2(62000-3420\pi+700e-860\kappa)+d(384384+15e^2\\
&&+15960e-30e\kappa-150\pi e+15\kappa^2+150\kappa\pi-75610\pi\\
&&-24340\kappa+375\pi^2)-400e^2-117120e+3920\pi e+960\kappa e\\
&&+226560\kappa-4720\kappa\pi-560\kappa^2+530880\pi-9600\pi^2.
\end{eqnarray*}

These calculations verify LeBarz' trisecant formula for $N_3$
\cite[Th\'eo\-r\`eme 8]{LeBarz} and the computation of $N_4$ by Tikhomirov and
Troshina \cite{TikhomirovN4}. The formula for $N_5$ seems to be new.
I omit the presentation of $N_6$ and $N_7$: the information is contained in
the following analysis of these numerical data. For $X=\IP^2$ and $\ko_X(H)=\ko_{\IP^2}(m)$ these tally with the polynomials computed by
Ellingsrud and Str\o{}mme using a torus action on $\IP^2$ and the Bott formula \cite{EllStrBott}.

Taking the logarithm of the generating function, we may write:
$$\sum_{n\geq0}N_nz^n=\exp\left(\sum_{m>0}\frac{(-1)^{m-1}}{m}d_mz^m\right)$$
where the coefficients $d_m$ {\sl a priori} are rational polynomials in $d=H^2$, $\pi=HK$, $\kappa=K^2$ and $e$. One can show that these polynomials are in
fact linear (cf.\ \cite{EllingsrudGoettscheLehn}). The explicit calculation yields
$$
\begin{array}{c}
d_1=d\phantom{-\pi-e}\\
d_2=10d+5\pi-e+\kappa\\
d_3=112d+96\pi-20e+28\kappa\\
d_4=1320d+1507\pi-324e+550\kappa\\
d_5=16016d+22120\pi-4880e+9440\kappa\\
d_6=198016d+314738\pi-70976e+151260\kappa\\
d_7=2480640d+4402720\pi-1012032e+2326192\kappa.\\
\end{array}
$$
From this one can attempt to guess the generating functions.
Let
$$k=z-9z^2+94z^3-\ldots\in\IQ[[z]]$$
be the inverse power series of the rational function
$$z=\frac{k(1-k)(1-2k)^4}{(1-6k+6k^2)^3}.$$
This is a solution of the differential equation
$$\frac{dz}{z}=\frac{dk}{k(1-k)(1-2k)(1-6k+6k^2)}.$$

\begin{conjecture}--- Using the notations above the following formula holds:
$$\sum_{n\geq0}N_nz^n=\frac{(1-k)^a\cdot(1-2k)^b}{(1-6k+6k^2)^{c}}$$
with $a={HK-2K^2}$, $b={(H-K)^2+3\chi(\ko_X)}$, and
$c=\frac{1}{2}H(H-K)+\chi(\ko_X)$.
\end{conjecture}

We thank Don Zagier for pointing out to us the existence of Sloane's `Encyclopedia of Integer Sequences' \cite{Sloane}. Intensive use of the 
on-line version of the Encyclopedia, numerous numerological experiments
and some inspiring help from Don Zagier allowed me to guess the 
generating functions. He also found a simple substitution to turn my still
awkward version of the generating function into the smooth form presented
above.


\subsection{The cohomology ring of $(\IA^2)^{[n]}$}


In this section we will describe an identification of the cohomology ring of
$(\IA^2)^{[n]}$ with the ring of certain explicitly given differential operators
on the polynomial ring in countably many variables. 

Of course, the affine plane $\IA^2$ is
not projective, so that we cannot directly apply the methods of the previous
sections. On the other hand, in \cite{Nakajima} Nakajima does work with 
non-projective surfaces, the only difference being that the operators $\qq_n$,
$n<0$, must be modelled on cohomology classes with compact support rather than
ordinary cohomology classes. The reason for this is that, in the notations of
Definition \ref{definitionofoscillators}, the morphism $p_1$ is proper, so that
push-forward is defined, whereas $p_2$ is proper only if the variety $X$ is
proper. With this modification Nakajima's main theorem holds for the affine plane as well.

As $H^*(\IA^2;\IQ)=\IQ$, we simplify notations by putting $q_m:=\qq_m(1_{\IA^2})$.
Then $\IH=\bigoplus_{n,i}H^i((\IA^2)^{[n]};\IQ)\isom \IQ[q_1,q_2,\ldots]$, the polynomial ring in countably infinitely many variables, and if $q_m$ is given 
degree $m$, then $\IH_n:=H^*((\IA^2)^{[n]};\IQ)$ is the homogeneous component of 
$\IH$ of degree $n$. As any vector bundle on $\IA^2$ is trivial, there is essentially only one tautological bundle $\ko^{[n]}$ on $(\IA^2)^{[n]}$. Let 
$\ch_i:\IH\to\IH$ be the components of the associated
Chern character operator, and let $\dd=\ch_1$ as before. The inclusion 
$\IA^2\subset\IP^2$ induces an open embedding $(\IA^2)^{[n]}\subset(\IP^2)^{[n]}$ which in turn gives rise to an epimorphism
of rings $H^*((\IP^2)^{[n]};\IQ)\to H^*((\IA^2)^{[n]};\IQ)$. This implies that
all commutation relations for the $q_m$ and $\ch_i$ hold in $\IH$ as well. In
fact they become much simpler as the pull-back both of $c_1(\IP^2)$ and $c_2(\IP^2)$ is zero.
To describe these relations in the given special setting, let $\partial_m:=m\frac{\partial}{\partial q_m}$. 

\begin{theorem}\label{chDidentification}--- The Chern character of the tautological bundle acts on $\IH$
as follows:
$$\ch_\nu=\frac{(-1)^\nu}{(\nu+1)!}\sum_{n_0,\ldots n_\nu>0}q_{n_0+\ldots+n_\nu}\partial_{n_0}\cdot\ldots\cdot\partial_{n_\nu}.$$
For each $n$, the cohomology ring $\IH_n$ is generated as a $\IQ$--algebra
by $ch_\nu(\ko^{[n]})$, and the relations between these generators are those
of the restriction of the given differential operators to $\IH_n$.
\end{theorem}

That $\IH_n$ is generated by the chern classes of the tautological bundle
had earlier been proved by Ellingsrud and Str\o{}mme \cite{EllStr2}.

In order to prove the theorem we consider a larger class of differential
operators on $\IH$ defined by
$$D_{n,\nu}:=\sum_{n_1,\ldots,n_\nu>0}q_{n+\sum_i n_i}
\prod_{i=1}^\nu\partial_{n_i}$$
for nonnegative integers $n$ and $\nu$, with the usual conventions
$D_{n,0}=q_n$ for $n>0$ and $D_{0,0}=0$.
The key observation is that $\dd=-\frac{1}{2}D_{0,2}$. This follows 
directly from Theorem \ref{maintheorem} and the fact that in the present
situation $H^*(\IA^2;\IQ)=\IQ$. It is easy to check by explicit calculation
that these operators satisfy the following commutation relations
$$[D_{n,\nu},D_{m,\mu}]=(\nu m-\mu n)\cdot D_{n+m,\nu+\mu-1}.$$
In particular, $q_n'=-\frac{1}{2}[D_{0,2},D_{n,0}]=-n\cdot D_{n,1}$, or more generally, by induction:
$$q_n^{(\nu)}=(-n)^{\nu}\cdot D_{n,\nu}.$$
We can now easily generalise Theorem \ref{cqcommutator}:
$$[\ch_n,q_m]=\frac{(-1)^n}{n!}m\cdot D_{m,n}.$$
For $m=1$ the assertion follows from the basic relation (Theorem \ref{cqcommutator}) $$[\ch_n,q_1]=\frac{1}{n!}q_1^{(n)}=\frac{(-1)^n}{n!}D_{1,n},$$
and for $m>1$ we deduce it by induction using
$-mq_{m+1}=[q_1',q_m]$ as well as $q_1'=-D_{1,1}$ and
$[\ch_n,q_1']=[\ch_n,q_1]'$.

\medskip

\prfofthetheorem We must first show that 
$$\ch_n=(-1)^n\frac{D_{0,n+1}}{(n+1)!}.$$
Observe that by the commutation rules for the operators $D_{*,*}$ we have
$$[\frac{(-1)^n}{(n+1)!}D_{0,n+1},q_m]=[\frac{(-1)^n}{(n+1)!}D_{0,n+1},D_{m,0}]
=\frac{(-1)^n}{n!}m\cdot D_{m,n}.$$
Thus $\ch_n$ and $\frac{(-1)^n}{(n+1)!}D_{0,n+1}$ show the same commutation
behaviour with all generators $q_m$ of $\IH$ and clearly act trivially on the
vacuum. Hence they are equal.

It remains to check that the Chern classes of the tautological bundle generate
$\IH_n$. Let $\lambda=(\lambda_1,\lambda_2,\ldots)$ be a partition of $n$,
i.e.\ $n=\|\lambda\|:=\sum_ii\lambda_i$ and let
$q_\lambda:=\prod_iq_i^{\lambda_i}$ be the associated monomial.
The monomials $q_\lambda$ with $\|\lambda\|=n$ form a $\IQ$-basis of $\IH_n$.
Let us say that $q_\lambda<q_\mu$ if $\lambda>\mu$ in the lexicographical order. We want to show that the subring $\IH_n'$ in $\IH$ generated
by the action of $\ch_m$, $m=1,\ldots,n-1$, on $1=q_1^n=q_{(n,0,\ldots)}$ 
contains the monomial $q_\lambda$ for all partitions $\lambda$ of $n$. As
this is true for the smallest possible monomial $q_{(n,0,\ldots)}$, we proceed
by induction. Given $\lambda$ we assume that $q_{\mu}\in\IH_n'$ for all $q_\mu<q_\lambda$. As $\lambda\neq(n,0,\ldots)$,
let $a$ be the smallest index $>1$ such that $\lambda_a>0$, i.e.\ $\lambda=(\lambda_1,0,\ldots,0,\lambda_a,\lambda_{a+1},\ldots).$
Consider now the partition
$$\lambda':=(\lambda_1+a,0,\ldots,0,\lambda_a-1,\lambda_{a+1},\ldots).$$
Then $q_{\lambda'}<q_{\lambda}$ and hence is contained in $\IH_n'$ and
$$ch_{a-1}q_{\lambda'}=(-1)^{a-1}\binom{\lambda_1+a}{a}q_\lambda+\ldots$$
where $\ldots$ stands for a linear combination of smaller monomials.
This finishes the induction.\qed

%
%

\noindent
Manfred Lehn\\
Mathematisches Institut der Georg-August-Universtit\"at\\
Bunsenstra\ss e 3-5, D-37073 G\"ottingen, Germany\\
e-mail: \verb{lehn@uni-math.gwdg.de{

 \end{document}